\documentclass[11 pt]{article}  
\usepackage{amssymb, latexsym}  
\usepackage{amsmath}  
\usepackage{mathrsfs}  
\begin{document}
  
\newtheorem{Lemma}{Lemma}[section]  
\newtheorem{Theorem}[Lemma]{MainTheorem}  
\newtheorem{theorem}[Lemma]{Theorem}  
\newtheorem{proposition}[Lemma]{Proposition}  
\newtheorem{corollary}[Lemma]{Corollary}  
\newtheorem{definition}[Lemma]{Definition}  
\newtheorem{remark}[Lemma]{Remark}  
\newtheorem{remarks}[Lemma]{Remarks}  
\newtheorem{Ex}[Lemma]{Example}  
\newtheorem{Exs}[Lemma]{Examples}  
\newtheorem{lemma}[Lemma]{Lemma}  
\newenvironment{Proof}{{\sc Proof.}\ }{\ \rule{1ex}{1ex}\vspace{0.5truecm}}
\newenvironment{SProof}{{\sc Sketch of the Proof.}\ }{\ \rule{1ex}{1ex}\vspace{0.5truecm}}  
\newcommand{\entrar}{\makebox[\parindent]{ }}  
\newtheorem{exemple}[Lemma]{Example}  

\newcommand{\Hom}{\mbox {\rm Hom}}  
\newcommand{\Spec}{\mbox{\rm Spec}}  
\renewcommand{\dim}{\mbox{\rm dim}}  
  
\title{{\bf Logarithmic  De Rham, Infinitesimal and Betti Cohomologies}}  
  
\author{Bruno Chiarellotto, Marianna Fornasiero,\\  
Dipartimento di Matematica, Universit\`a degli Studi di Padova,\\  
Via Belzoni 7, 35100 Padova, Italy \\  
email: chiarbru@math.unipd.it,  mfornas@math.unipd.it \thanks{
{\sl Key words}: 
Log Schemes,  Log  De Rham Cohomology, Log Betti Cohomology,  Log Infinitesimal Cohomology. {\sl Supported by}: EC, RTN (Research Training Network) Arithmetic Algebraic Geometry (AAG); MIUR, GVA project; Universit\`a di Padova PGR ``CPDG021784". {\sl 2000 AMS Classification}: 14FXX.}}  
  
\maketitle


\section*{Introduction}

For a non singular scheme $Y$ over $\mathbb C$ the hyper-cohomology of the
algebraic de Rham complex calculates the analytic cohomology
$H^{^.}(Y^{an}, \mathbb C)$. 
  For singular $Y$ there is no straighforward
generalization of this calculation:  indeed,
it is  the algebraic side that causes problems.  
Grothendieck has  introduced the algebraic  infinitesimal
site  $Y_{inf}$  \cite{Groth1}.  Moreover,  as explained in  \cite{Hartshorne}, when $Y$ admits an  embedding
 as a closed subscheme of a
smooth scheme $X$, one can also consider the completion of the de Rham complex
$\Omega^{^.}_ {X}$ along $Y$: $\Omega^{^.}_
{X\hat{|}Y}$. At this point one has three different cohomologies 

\begin{equation}
\label{01}
\mathbb H^{^.}(Y, \Omega^{^.}_ {X\hat{|}Y}),
\end{equation}

\begin{equation}
\label{02}
\mathbb
 H^{^.}(Y_{inf}, \mathcal O_{Y_{inf}}),
\end{equation}
and
\begin{equation}
\label{03}
H^{^.}(Y^{an}, \mathbb C).
\end{equation}

The  isomorphism between (2) and (3) was proved by Grothendieck
(\cite{Groth1}) only in the case of a
smooth scheme over $\mathbb C$.  The
  isomorphism between (1) and (3)  was proved by Herrera-Liebermann
(\cite{Herrera}), in the case of $Y$ proper over
$\mathbb C$, while  Deligne (unpublished), and Hartshorne (\cite[Chapter IV,
Theorem (I.I)]{Hartshorne}) proved it for a general (not
necessary proper) scheme over
$\mathbb C$.  A direct statement asserting 
the isomorphism of these cohomology groups for arbitrary $\mathbb C$-schemes $Y$
 cannot be found in literature, although all the necessary ingredients are
given.  The proof presented in this paper, if
applied to classical schemes, can be used to fill this gap (see \S 1).  Of course the generalization of this problem to the case of mixed or
finite characteristic  has been carefully studied by Berthelot and Ogus.\\

On the other hand, in more recent years the notion of scheme and the properties
of schemes  have been generalized by the introduction of log schemes. Among the expected features of log schemes,
 there is the fact that log smooth  schemes
(which are in general singular as schemes) should behave like classical
smooth schemes and moreover should also be related
 to analytic schemes.  The goal of the
present work is to introduce the
log scheme analogues of (1),(2),(3)  over $\mathbb C$,  and  prove the isomorphisms between them.

With these ideas in mind we  first consider the analogue of Grothendieck's Infinitesimal
Site (\cite{Groth1}) in the logarithmic context (see also \cite{Kato1} for positive characteristic). 
We work with pro-crystals and we link them to the logarithmic stratification on 
pro-objects. If we suppose that there exists an exact closed
immersion of fs log schemes $Y \hookrightarrow X$, with $X$
log smooth over $\mathbb C$, then  we can define the Log De Rham Cohomology of $Y$ over $\mathbb C$, as the
hyper-cohomology of the complex $\omega^{^.}_{X \hat |
Y}$, where $X \hat | Y$ represents the formal completion of the log scheme
$X$ along its closed log subscheme $Y$. We give a
direct proof of the {\sl existence of an isomorphism between
the Log Infinitesimal Cohomology of $Y$ over $\mathbb C$, and its
Log De Rham Cohomology}, namely we prove the following isomorphism
$$
H^{^.}(Y^{log}_{inf}, \mathcal O_{Y^{log}_{inf}}) \cong
H^{^.}_{DR,log}(Y/\mathbb C)=: \mathbb H^{^.}(Y, \omega^{^.}_{X \hat{|}Y})
$$
This result was proved by Shiho (\cite{Shiho}) for a log smooth log scheme
over $\mathbb C$. \\

For the remaining  isomorphisms,  we were inspired by a work of K. Kato and
C. Nakayama
(\cite[Theorem (0.2), (2)]{Nakayama}).   Given  an fs
(ideally) log  smooth  log scheme $X$  over $\mathbb C$, Kato and Nakayama associate
a   topological space $X^{an}_{log}$
 and show that the
algebraic Log De Rham cohomology of $X$
(which is
defined as
the hyper-cohomology of the log De Rham complex $\omega^{^.}_X$) is
isomorphic to the cohomology of the constant sheaf $\mathbb
C$ on $X^{an}_{log}$, i.e.
\begin{equation}
\label{0}
H^{^.}_{DR,log}(X/\mathbb C)=: \mathbb H^{^.}(X, \omega^{^.}_X) \cong
H^{^.}(X^{an}_{log}, \mathbb C)
\end{equation}
We will prove an analogue of \eqref{0},  for a general  fs log scheme
$Y$ over $\mathbb C$.  In this case, the Log De Rham Cohomology of
$Y$ is defined  as before (for $Y$ admitting an exact closed immersion).
In \S \ref{o}, we
extend the theory of Kato-Nakayama (\cite{Nakayama}) to
the log formal setting.   To this end, we first introduce a ringed
topological space $(X \hat{|}Y)^{log}$, associated to the log
formal analytic space $(X \hat{|}Y)^{an}$, with sheaf of rings $\mathcal
O^{log}_{(X \hat{|}Y)^{an}}$ (Definition \ref{re4}).
The underlying topological space $(X \hat{|}Y)^{log}$ of this ringed space
coincides with $Y^{an}_{log}$. We construct the
complex $\omega^{^.,log}_{(X \hat{|}Y)^{an}}$ (Definition \ref{re9} and
\eqref{re11}), which is a sort of ``formal analogue" of
the complex $\omega^{^.,log}_{X^{an} }$, introduced  by Kato-Nakayama for a
log smooth log scheme $X$ (\cite[(3.5)]{Nakayama}).

Later, in \S \ref{formalres}, we give a ``formal version" of the Deligne
Poincar\'e Residue map (\cite[(3.6.7.1)]{Deligne2}), in
the particular case of a smooth scheme $X$ over $ \mathbb C$,
endowed with log structure given by a normal crossing
divisor $D$, and $Y \hookrightarrow X$ a closed subscheme,  with the
induced log structure (\S \ref{res9}).  We show that this
map is an isomorphism. It is useful for describing the cohomology of the
complex $\omega^{^.}_{(X \hat{|}Y)^{an}}$.\\  Using that 
description, we can prove the  Log Formal Poincar\'e Lemma (Theorem
\ref{poin0}):  {\sl given a general fs log scheme $Y$ over
$\mathbb C$, the Betti Cohomology of its associated topological space
$Y^{an}_{log}$ is isomorphic to the hyper-cohomology of
the complex $\omega^{^.,log}_{(X \hat{|}Y)^{an}}$}.

 In \S \ref{finale}, we
show that $\omega^{^.}_{(X \hat{|}Y)^{an}}$ is
quasi-isomorphic to $\mathbb R \tau_* \omega^{^.,log}_{(X \hat{|}Y)^{an}}
\cong \mathbb R \tau_* \mathbb C_{Y^{log}}$
(Proposition \ref{f0}), where $\tau \colon Y^{an}_{log} \longrightarrow
Y^{an}$ is the canonical (continuous, proper and
surjective) Kato-Nakayama map of topological spaces. Finally, we prove that
there exists an  isomorphism in cohomology $\mathbb
H^{^.}(Y^{an}, \omega^{^.}_{(X \hat{|}Y)^{an}}) $ \\ $\cong \mathbb H^{^.}( Y,
\omega^{^.}_{X \hat{|}Y})$ between the Analytic and the
Algebraic Log  De Rham Cohomology (Theorem  \ref{f3}).  We conclude with
the main theorem of this article (Theorem \ref{f00}):
{\sl the cohomology of the constant sheaf $\mathbb C$ on the topological
space $Y^{an}_{log}$, associated to an fs log scheme
$Y$, is isomorphic to the Log  De Rham Cohomology of $Y$}, 
$$
H^{^.}(Y^{log}_{inf}, \mathcal O_{Y^{log}_{inf}}) \cong H^{^.}_{DR,log}(Y/\mathbb C) 
\cong H^{^.}(Y^{an}_{log}, \mathbb C)
$$

\vskip20pt

We would like to thank Prof. L. Illusie and T. Tsuji,  for their   precious comments and suggestions.

\vskip20pt

\section*{Preliminaries} 

{\bf Notations}: by  $S$ we denote the logarithmic scheme $\Spec \, \mathbb C$ endowed with the trivial log structure, and, by a {\sl log scheme}, we mean a logarithmic scheme over $S$, whose underlying scheme is a separated $\mathbb C$-scheme of finite type. 
Moreover, if $A^{^.} $ is a complex of sheaves and $k \in \mathbb N$, then $A^{^.}[-k]$ is the complex defined in degree $j$ as $A^{j+k}$. 
  
\subsection{The Logarithmic Infinitesimal Site}
\label{sito}

Given a log scheme  $X$, endowed with a fine log structure $M$, we denote by $InfLog(X/S)$ the Logarithmic Infinitesimal Site of $X$ over $S$. It is given by $4-$uples $(U,T,M_T,i)$, where $U$ is an \'etale scheme over $X$, $(T,M_T)$ is a scheme with a fine log structure   
over $S$, $i$ is an exact closed immersion $(U,M) \hookrightarrow (T,M_T)$ over $S$, defined by a  
nilpotent ideal on $T$, i.e. $i$ is a nilpotent exact closed immersion. Morphisms, coverings (for the usual \'etale topology), and sheaves on $InfLog(X/S)$ are defined  in the usual way. The category of all sheaves on  
$InfLog(X/S)$ is a ringed topos,  called the Logarithmic  Infinitesimal Topos  
of $X$ over $S$, and denoted by $(X/S)^{log}_{inf}$.  

\subsection{Pro-Crystals and Logarithmic Stratification}
\label{stratcrys0}

Let $X$ be a log smooth log scheme. For the definition of  pro objects we refer to \cite{Bald}, \cite{Fiorot}, \cite[\S 6.2]{Groth1}. 
\begin{definition}  
\label{crystal000}
Let $\{  \mathscr F_k \}_{k \in K}$ be a pro object. It is a {\bf  pro-crystal}  if, given $g \colon (U^{\prime}, T^{\prime}, M_{T^{\prime}}, i^{\prime}) \longrightarrow (U,T,M_T,i)$ a morphism in $X^{log}_{inf}$,  there exists an isomorphism of pro objects  
$$
\{ g^* \mathscr F_k (U,T,M_T,i) \}_{k \in K} \cong \{  \mathscr F_k (U^{\prime}, T^{\prime}, M_{T^{\prime}}, i^{\prime}) \}_{k \in K} 
$$
i.e. $\{ g^* \mathscr F_{kT}  \}_{k \in K} \cong \{  \mathscr F_{k T^{\prime}} \}_{k \in K} $. In a similar way we can define  Artin-Rees pro-crystals (see \cite[Proposition 0.5.1]{Ogus}).
\end{definition}

For each integer $i \geq 0$, let $\Delta^1_{log}(i)$ be the $i-$th log infinitesimal neighbourhood of the diagonal $(X,M) \hookrightarrow (X,M) \times_{S} (X,M)$, and let $\Delta^2_{log}(i)$ be the $i-$th log infinitesimal neighbourhood of $(X,M) \times_{S} (X,M) \hookrightarrow (X,M) \times_{S} (X,M) \times_{S} (X,M) $ (where the fiber   product is taken in the category of fine log schemes). We have the canonical projections $p_1(i), p_2(i) \colon \Delta^1_{log}(i) \longrightarrow (X,M) $, and $p_{31}(i), p_{32}(i), p_{21}(i) \colon \Delta^2_{log}(i) \longrightarrow \Delta^1_{log}(i)$. We denote by $\mathscr P^{\nu,i}_{X,log}$ the structural sheaf of rings $\mathcal O_{\Delta^{\nu}_{log}(i) }$, for each $\nu =1,2$, $i \geq 0$. In particular, we can regard $\mathscr P^{1,i}_{X,log} $ as an $\mathcal O_X-$module in two ways, via  the canonical projections $p_1(i), p_2(i)$. So, we call the left $\mathcal O_X-$module structure (resp. right $\mathcal O_X-$module structure) on $\mathscr P^{1,i}_{X,log} $ the structure given by $p_1(i)$ (resp. $p_2(i)$). \\

We introduce a logarithmic stratification  on the category of pro-coherent $\mathcal O_X-$modules. We could define a logarithmic stratification 
``at any level" of the pro object, and consider the  pro-category of log stratified $\mathcal O_X-$modules. But this stratification would be too
 restrictive for our  purpose. We need to work with a larger category and, to this aim, we introduce the logarithmic stratification  as a pro-morphism. 
\begin{definition}
\label{strat01} \cite[Definition 1.3]{Fiorot}
Let $\{ \mathscr F_k \}_ {k \in K}$ be a pro-coherent  $\mathcal O_X$ module. A {\bf logarithmic stratification} on $\{ \mathscr F_k \}_ {k \in K}$ is a pro-morphism 
$$
\{ \mathscr F_k \}_{k} \stackrel{s_{\{ \mathscr F_k \}_{k}}}{ \longrightarrow} \{    \mathscr F_k \}_{k} \otimes \{ \mathscr P^{1,i}_{X,log} \}_i 
$$
such that the coidentity diagram 
$$
\begin{matrix}
\{ \mathscr F_k \}_ {k } &  \stackrel{s_{\{ \mathscr F_k \}_{k}}}{ \longrightarrow} & \{    \mathscr F_k \}_{k} \otimes \{ \mathscr P^{1,i}_{X,log} \}_i \\
^{id_{\{    \mathscr F_k \}_{k}}  }{\downarrow} & & {\downarrow}^{id_{\{    \mathscr F_k \}_{k}} \otimes \{ q_{i,0 } \}_i} \\
\{ \mathscr F_k \}_ {k }  & \stackrel{id_{\{ \mathscr F_k \}_ {k }}}{\longrightarrow} & \{ \mathscr F_k \}_ {k } 
\end{matrix}
$$
and the coassociativity diagram 
$$
\begin{matrix}
\{ \mathscr F_k \}_ {k } &  \stackrel{s_{\{ \mathscr F_k \}_{k}}}{ \longrightarrow} & \{    \mathscr F_k \}_{k} \otimes \{ \mathscr P^{1,i}_{X,log}
 \}_i \\
^{s_{\{    \mathscr F_k \}_{k}}  }{\downarrow} & & {\downarrow}^{s_{\{    \mathscr F_k \}_{k}} \otimes id_{\{ \mathscr P^{1,i}_{X,log} \}_i}}\\
\{ \mathscr F_k \}_ {k }  \otimes \{ \mathscr P^{1,i}_{X,log} \}_i  & \stackrel{id_{\{ \mathscr F_k \}_ {k }} \otimes s_{\{ \mathscr P^{1,i}_{X,log}
 \}_i}}{\longrightarrow} & \{ \mathscr F_k \}_ {k }  \otimes \{ \mathscr P^{1,i}_{X,log} \}_i  \otimes \{ \mathscr P^{1,i}_{X,log} \}_i 
\end{matrix}
$$
are commutative, where $ q_{i,j } \colon \mathscr P^{1,i}_{X,log} \longrightarrow \mathscr P^{1,j}_{X,log}$ are the natural compatible maps, 
and $s_{\{ \mathscr P^{1,i}_{X,log} \}_i} = \{ \delta_X^{i,j} \}_{(i,j)}$ (see \cite[Lemma 3.2.3]{Shiho} for the definition of  $\delta_X^{i,j} \colon \Delta^1_{log}(i) \times_{(X,M)}  \Delta^1_{log}(j) \longrightarrow  \Delta^1_{log}(i+j)$). 
\end{definition}
\begin{theorem}
\label{main0} \cite[Theorem 3.2]{Fiorot}. 
The following two categories are equivalent:\\
(a) the category of pro-crystals on $InfLog(X/S)$;\\
(b) the category of  pro-sheaves $\{ \mathscr M_n\}_n$ on $X$, endowed
with a logarithmic stratification.
\end{theorem}
\begin{remark}
\label{strat03}
In fact, our pro-crystals  are actually Artin-Rees pro-crystals and one could refine the previous result on these objects. 
\end{remark}
Now let $i \colon Y \hookrightarrow X$ be an exact closed immersion of fs log schemes, with  $X$ log smooth over $S$.   We consider the direct image 
functor $i^{log}_{inf*} \colon Y^{log}_{inf} \longrightarrow X^{log}_{inf}$. 
For a crystal $\mathscr F$ of $Y^{log}_{inf}$, we briefly describe  the construction of the direct image $i^{log}_{inf*} \mathscr F$, 
in characteristic zero. Let $(U,T,M_T,j) \in InfLog((X,M)/S)$, then we consider the fiber product (in the category of fine log schemes) 
$U_Y= (Y,N) \times_{(X,M)} (U,M)$. By base change, the map $U_Y \hookrightarrow (U,M)$ is an exact  closed immersion. 
Since $j \colon (U,M) \hookrightarrow (T,M_T)$ is a nilpotent exact closed immersion, $U_Y \hookrightarrow (T,M_T)$ 
is also an exact closed immersion, and we can take the $n-$th log infinitesimal neighbourhood of $U_Y$ inside $(T,M_T)$, 
and denote it by $(T_n,M_n)$. Let $\lambda_n \colon (T_n,M_n) \hookrightarrow (T,M_T)$. Then, 
$(i^{log}_{inf*} \mathscr F)_{(U,T,M_T,j)}=:  \displaystyle \lim_{\stackrel{\longleftarrow}{n}}  \lambda_{n*}  \mathscr F_{(U_Y,T_n,M_n, j_n)} $. 
So, the Artin-Rees  pro-crystal $\{  \mathscr F_n \}_{n \in \mathbb N}$ on $InfLog((X,M)/S)$,  
associated to $i^{log}_{inf*} \mathscr F$,  is defined on $(U,T,M_T,j) \in InfLog((X,M)/S)$ as 
$\mathscr F_{n (U,T,M_T,j)}:= \lambda_{n*} (\mathscr F_{(U_Y,T_n,M_n, j_n)} )= \lambda_{n*} (\mathscr F_{T_n})$ (\cite[Proposition 0.5.1]{Ogus}). 
In particular, the Artin-Rees pro-crystal $\{ \mathcal O_{n} \}_{n \in \mathbb N}$ associated to  $i^{log}_{inf*} \mathcal O_{Y^{log}_{inf}} $, 
 is in fact defined, on each $(U,T,M_T,j)$,  as 
$( \{  \mathcal O_{n} \}_{n \in \mathbb N} )_{(U,T,M_T,j)}:=\{  \lambda_{n*} \mathcal O_{T_n} \}_{n \in \mathbb N}$. 
\begin{remark}
\label{re00}
{\rm From Theorem \ref{main0}, the log stratified $\mathcal O_X-$pro-module associated to the (Artin-Rees)  pro-crystal 
$\{ \mathcal O_n\}_{n \in \mathbb N} $ is equal to $\{ (\mathcal O_n)_{(X,X,M,id)}\}_{n \in \mathbb N} = 
\{ \lambda_{n*} \mathcal O_{Y_n} \}_{n \in \mathbb N}  = \{ \mathcal O_{X}/\mathscr I^n\}_{n \in \mathbb N} $, 
where  $\mathscr I \subset \mathcal O_X$ is the ideal of definition of $Y$ into $X$, and $Y_n$ is the $n-$th 
log infinitesimal neighborhood of $Y \hookrightarrow X$. Moreover, $(i^{log}_{inf*} \mathcal O_{Y^{log}_{inf}})_{(X,X,M,id)} 
= \displaystyle \lim_{\stackrel{\longleftarrow}{n}}  \mathcal O_{X}/\mathscr I^n =  \mathcal O_{X \hat{|}Y}$. }
\end{remark}

\subsection{Linearization of the Log De Rham complex}
\label{linearizzato0}

Let $\omega^{^.}_X$ be the log De Rham complex of the log smooth log scheme $X$ over $S$. 
As in the classical case (\cite[p. 2.17]{Berth-Ogus}), we denote the complex of Artin-Rees pro-coherent $\mathcal O_X$ modules which is the linearization of $\omega^{^.}_X$  by $\{ L_X(\omega^{^.}_X)_i \}_{i \in \mathbb N}$,  i.e. 
\begin{equation}
\label{eq1}
\{ L_X(\omega^{^.}_X)_i \}_{i \in \mathbb N} =: \{ \mathscr P^{1,i}_{X,log} \}_{i \in \mathbb N}  \otimes_{\mathcal O_X} \omega^{^.}_X
\end{equation}
Since, for all $i,j \in \mathbb N$, there exist maps $\mathscr P^{1,i+j}_{X,log} \otimes \omega^{k}_X \longrightarrow \mathscr P^{1,i}_{X,log} \otimes  \omega^{k}_X \otimes \mathscr P^{1,j}_{X,log}$, each  term of \eqref{eq1} has a canonical logarithmic stratification, in the sense of Definition \ref{strat01} (\cite[Construction 2.14]{Berth-Ogus}).\\ We have a  local description of the differential maps $ \{ L_X(d^{^.}_X)_i \}_i $ of this complex. Indeed, let  $M^i$ be the log structure on $\Delta^{1}_{log}(i)$. Let $U \longrightarrow X$ be an \'etale morphism of schemes, and let $m \in \Gamma(U, M)$. Then, there exists uniquely an element $u_{m,i} $ in $\Gamma(U, (\mathscr P^{1,i}_{X,log})^* ) \subset \Gamma (U, M^i)$, such 
that $p_2(i)^* (m)= p_1(i)^*(m) u_{m,i}$ (\cite[pp. 43, 44]{Shiho}). 
In particular, we have that $u_{m,i} -1 \in {\rm Ker} \, \{ \Gamma(U, (\mathscr P^{1,i}_{X,log})^* )  \longrightarrow \Gamma (U, \mathcal O_X) \}$ 
(\cite[Lemma 3.2.7]{Shiho}).  \\
Let now $x \in X$, and  $t_1,...,t_r \in M_{x}$ be such that $\{ {\rm dlog} \, t_j \}_{1 \leq j \leq r }$ 
is a basis of $\omega^1_{X,x}$. We can restrict to an \'etale neighborhood $U$ of $x$, and suppose that $\{ {\rm dlog} \, t_j \}$
 is a local basis of $\omega^1_{X}$ on $U$. Let $u_{j,i}$, $1 \leq j \leq r$, $i \geq 0$, be the elements in $\Gamma(U, (\mathscr P^{1,i}_{X,log})^* )$ 
such that  $p_2(i)^* (t_j)= p_1(i)^*(t_j) u_{j,i}$, as above. We put $\xi_{j,i} := u_{j,i} -1 \in \Gamma(\Delta^1_{log}(i), \mathscr P^{1,i}_{X,log} )$
 (note that the $\xi_{j,i}$'s are compatible with respect to $i$). 
\begin{proposition}
\label{c01}
In the above notations, \\
(1) \cite[Lemma 3.2.7]{Shiho}, for each $i \geq 0$, the $\mathcal O_{X}-$module $\mathscr P^{1,i}_{X,log}$ is locally free with basis 
$$
\{ \xi^a_i := \prod_{j=1}^r \xi^{a_j}_{j,i} | 0 \leq \sum_{j=1}^r a_j \leq i\}
$$
where $a=(a_1,...,a_r) $ is a multi-index of length $r$. In particular, $\{ \xi_{1,1},...,\xi_{r,1} \}$ is a basis for the 
locally free $\mathcal O_{X}-$module $\mathscr P^{1,1}_{X, log}$, \'etale locally at $x$;\\
(2) \cite[Proposition 3.2.5]{Shiho}, there exists a canonical isomorphism of $\mathcal O_X-$modules 
$ \mathscr K_1 / \mathscr K_1^2 \stackrel{\cong}{\longrightarrow} \omega^1_X $, 
where $ \mathscr K_1:= {\rm Ker} \, \{ \Delta^* \colon \mathscr P^{1,1}_{X,log} \longrightarrow \mathcal O_X \}$ 
(with $\Delta \colon X \hookrightarrow \Delta^1_{log} (1)$ the exact closed immersion). Under this identification,  
the local basis $\{ {\rm dlog} \, t_j \}_{1\leq j \leq r}$ of $\omega^1_X$ is identified with $\{ \xi_{j,1}  \}_{1\leq j \leq r}$. 
\end{proposition}
Therefore, for each $i \geq 0$, the map $ L_X(d^{^.}_X)_i $ is the $\mathcal O_X-$linear map 
$$
 L_X(d^{k}_{X})_i \colon \mathscr P^{1,i}_{X, log}    \otimes_{\mathcal O_X} \omega^k_{X} \longrightarrow  \mathscr P^{1,i-1}_{X, log}
    \otimes_{\mathcal O_X} \omega^{k+1}_{X}
$$ 
defined, for $a \in \mathcal O_{X}$, $\omega \in \omega^k_{X}$, and $n_j \in \mathbb N$ such that $n_1+...+n_r \leq i$, by setting
\begin{equation}
\label{equaz0}
L_X(d^{k}_{X})_i (a \xi_{1,i}^{n_1} \cdot ... \cdot \xi_{r,i}^{n_r}  \otimes \omega )=a \cdot \sum_{j=1,...,r} n_j \xi_{1,i-1}^{n_1}
 \cdot ... \cdot \xi_{j,i-1}^{n_j-1} \cdot ... \cdot \xi_{r,i-1}^{n_r} \otimes {\rm dlog} \, t_j  \wedge \omega + a \,  \xi_{1,i-1}^{n_1} 
\cdot ... \cdot \xi_{r,i-1}^{n_r} \otimes d \omega 
\end{equation}

\subsection{Log Formal Completion and Kato-Nakayama Topological Space}
\label{Weight0}  

Let  $X$ be a log scheme, endowed with fs log structure $M$, and let $i \colon Y \hookrightarrow X$ be an exact closed immersion. 
We refer to \cite[(5.8)]{Kato1} and \cite[Proposition-Definition 3.2.1]{Shiho}, for the definition of logarithmic $n-$infinitesimal 
neighborhoods ($n \in \mathbb N$). 
\begin{definition}
\label{def1}
The log formal completion $X \hat{|}Y$ of
$X$ along $Y$ is the classical formal completion of the scheme
$X$ along its closed subscheme $Y$, endowed with log structure given by the inverse image of $M$ via the canonical map
$X \hat{|}Y \longrightarrow X$.
\end{definition}
We note that, if the  closed immersion
$i \colon Y \hookrightarrow X$ is not exact, it is also possible to 
define the log formal completion of $X$ along $Y$,   
but in this article we will always work with exact closed immersions.\\
Let now $X^{an}$ be the (fs) log analytic space associated to $X$. Kato-Nakayama define the topological space $X^{an}_{log}$ associated 
to $X^{an}$ as the set 
$\{ (x,h) | x \in X^{an}, h \in {\rm Hom}(M^{gp}_x, \mathbb S^1), h(f)=f(x)/|f(x)|, $  
for any $f \in \mathcal O^*_{X^{an},x} \}$ (where $\mathbb S^1= \{ x \in \mathbb C; |x|=1 \}$). 
Let now $\beta \colon P \longrightarrow M$ be a fixed local chart for $X^{an}$, with $P$ an fs monoid.  
The topology on $X^{an}_{log}$ is locally defined as follows:
\begin{definition}
\label{logantopology}
In the local chart $\beta$,  $X^{an}_{log}$ is identified with a closed subset of 
$X^{an} \times {\rm Hom} (P^{gp}, \mathbb S^1)$,  
via the map $X^{an}_{log} \hookrightarrow X^{an} \times {\rm Hom}(P^{gp}, \mathbb S^1)$: $(x,h) \longmapsto (x, h_P)$,   
where $h_P$ is the composite $P^{gp} \longrightarrow M^{gp}_x \stackrel{h}{\longrightarrow } \mathbb S^1$.  
So, $X^{an}_{log}$ is locally endowed with the topology  induced from the natural topology on   
$X^{an} \times {\rm Hom}(P^{gp}, \mathbb S^1)$. 
\end{definition}
This local topology does not depend on the choice of the chart, so it induces a well defined global topology on $X^{an}_{log}$ 
(\cite[(1.2.1), (1.2.2)]{Nakayama}).
There exists a surjective map of topological spaces  $\tau \colon X^{an}_{log} \longrightarrow X^{an} \colon  (x,h) \longmapsto x$, 
which is continuous and proper 
(\cite[Lemma (1.3)]{Nakayama}).   
Though $X_{log}^{an}$ in general is not an analytic space, it is still endowed with a nice sheaf of rings
$\mathcal O^{log}_{X^{an}}$. Indeed, let $\mathcal L_X$ be  the sheaf of abelian groups on $X^{an}_{log}$ which 
represents the ``logarithms of local sections of $\tau^{-1}(M^{gp})$"  (\cite[(1.4)]{Nakayama}). 
There exists an exact sequence of sheaves of abelian groups
\begin{center}
$
\begin{array}{ccccccccc}
0 & \longrightarrow  & \tau^{-1}(\mathcal O_{X^{an}}) &  \stackrel{k}{\longrightarrow }
& \mathcal L_X & \stackrel{{\rm exp}}{\longrightarrow } & \tau^{-1}(M^{gp}/ \mathcal O^*_{X^{an}})
& \longrightarrow & 0
\end{array}
$
\end{center}
If we consider commutative $\tau^{-1}(\mathcal O_{X^{an}})$-algebras $\mathcal A$ on $X^{an}_{log}$, endowed with
a homomorphism $\mathcal L_X \longrightarrow \mathcal A$ of sheaves of abelian groups which commutes with $k$,  then 
$\mathcal O^{log}_{X^{an}}$ is the universal object among such $\mathcal A$ (\cite[(3.2)]{Nakayama}).\\
We suppose that $X$ satisfies the following hypothesis (\cite[Theorem (0.2), (2)]{Nakayama})

(*) Locally for the \'etale topology, there exists an fs monoid $P$, an ideal $\Phi$ of $P$, and a morphism 
$f \colon X \longrightarrow \Spec \, (\mathbb  C[P]/(\Phi))$ of log schemes over $S$, such that the underlying 
morphism of schemes is smooth,  and the log  
structure on $X$ is associated to $P \longrightarrow \mathcal O_X$.  
\begin{remark}  
\label{condition} {\rm We note that, if $X$ is (ideally) log smooth over  
$S$, then it satisfies hypothesis (*),   
because $X$ is a filtered semi-toroidal variety (\cite[Definition 5.2]{Ishida}, \cite[Definition (1.5)]{Illusie1}, \cite[Proposition II.1.0.11]{Marianna}).}  
\end{remark}  
\begin{theorem}  
\label{Comparison} \cite[Theorem (0.2), (2)]{Nakayama}.   
Let $X$ be an fs (ideally) log smooth log scheme over $S$ (see Remark \ref{condition}). Then, there exists a canonical isomorphism  
\begin{center}  
$  
\mathbb H^q(X, \omega^{^.}_X) \cong H^q(X^{an}_{log},\mathbb C), \text{ for all } q \in \mathbb Z.  
$  
\end{center}  
\end{theorem}

\section{Log Infinitesimal and Log De Rham Cohomologies}
\label{oo}

From now on, let $Y$ be  a generic fs log scheme over $S$,  endowed with log structure $M_Y$. We suppose there exists an exact closed 
immersion $i \colon Y \hookrightarrow X$ (i.e. $i^* M_X \cong M_Y$), where $X$ is an fs log smooth log scheme over $S$, with log structure $M_X$. 
The Log De Rham complex of the log
scheme $Y$ is the complex $(\omega^{^.}_{X \hat{|} Y} \cong \omega^{^.}_{X} 
\otimes_{\mathcal O_X} \mathcal O_{X \hat{|}Y}, \hat{d}^{^.})$, where $\hat{d}^{^.}$ is the
integrable connection induced by the differential $d^{^.}_X$ of $\omega^{^.}_{X} $. Then, we define the Log De Rham
Cohomology of $Y$ as the hyper-cohomology of
the log De Rham complex $\omega^{^.}_{X \hat | Y}$, i.e. $ H^{^.}_{DR,log}(Y/S)=: \mathbb H^{^.}(Y,  
\omega^{^.}_{X \hat | Y}) $. \\
Let $\mathscr I$ be the ideal of definition of $Y$ in $X$. For each fixed $\nu \geq 0$, we consider the diagonal 
immersion of fine log schemes $X \hookrightarrow X^{\nu}$, where $X^{\nu}$ is the fiber product over $S$ of $\nu +1 $ copies of
 $(X,M_X)$ over $S$ (the fiber product being that of the category of fine log schemes over $S$). We denote by $\Delta^{\nu}_{X,log}(i)$ the
 $i$-th log infinitesimal neighbourhood of the diagonal of $X^{\nu}$, and by $\mathscr P^{\nu,i}_{X,log}$ its structural scheaf of rings
 $\mathcal O_{\Delta^{\nu}_{X,log}(i)}= \mathcal O_{X^{\nu}} /\mathscr K_{\nu}^{i+1}$, where $\mathscr K_{\nu}$ is the ideal of
 definition of the log scheme $X$ inside $X^{\nu}$. \\
Now, if we fix $\nu$ and vary $i \in \mathbb N$, we get the  Artin-Rees  pro-object  of sheaves $\{ \mathscr P^{\nu,i}_{X,log} \}_i$ on $X$.  
On the other hand, if we fix $i$ and vary $\nu \in \mathbb N$, we get a sheaf on the simplicial log smooth log scheme 
$$
... \longrightarrow X^{\nu} 
\longrightarrow ... \longrightarrow 
X^1=X \times_S X \longrightarrow X
$$
which is the following cosimplicial sheaf of rings on $X $
\begin{equation}
\label{o2}
0 \longrightarrow  \mathcal O_{X} \stackrel{d_1-d_0}{\longrightarrow} \mathscr P^{1,i}_{X,log}  \stackrel{d_2-d_1+d_0}{\longrightarrow} \mathscr P^{2,i}_{X,log} \longrightarrow ... \longrightarrow \mathscr P^{\nu,i}_{X,log}  \longrightarrow ...
\end{equation}
where the maps are given by the alternating sum of the faces of the simplicial log scheme $\{ X^{\nu} \}_{\nu}$. If we vary $\nu$ and $i$,  we get a cosimplicial sheaf of Artin-Rees $\mathcal O_X-$pro modules $\{ \mathscr P^{\nu,i}_{X,log}  \}_{\nu,i}$. \\
We define the cosimplicial Artin-Rees pro-object $\{ Q_{log}^{\nu,i} \}_{\nu,i}$, by setting
\begin{equation}
\label{o3}
Q_{log}^{\nu,i}=: \mathscr P_{X, log}^{\nu+1,i}
\end{equation}
for every $i, \nu \geq 0$. Then, for each $\nu \geq 0$, there is a canonical homomorphism of pro-rings $\alpha_{log}^{\nu,^.} \colon \mathscr  P_{X, log}^{\nu,^.} \longrightarrow Q_{log}^{\nu,^.}$, defined by the canonical injection $\{ 0,1,...,\nu \} \hookrightarrow \{ 0,1, ..., \nu, \nu+1\}$. So, we have a homomorphism of cosimplicial pro-rings
\begin{equation}
\label{o4}
 \{ \alpha_{log}^{*,i} \}_i \colon \{ \mathscr P_{X,log}^{*,i} \}_i  \longrightarrow  \{ Q_{log}^{*,i} \}_i
\end{equation}
Let $\mathscr N$ be an $\mathcal O_X-$module. As in the classical case, we define the cosimplicial pro-ring 
\begin{equation}
\label{o5}
\{ Q_{log}^{*,i}(\mathscr N) \}_i=: \{ Q_{log}^{*,i} \}_i \otimes_{\mathcal O_X}  \mathscr N 
\end{equation}
We note that, for fixed $\nu \geq 0$, $\{ Q_{log}^{\nu, i}( \mathscr N) \}_i$ is clearly a $\{ Q_{log}^{\nu,i} \}_i-$module, and so an $\mathcal O_X-$bimodule with the obvious left and right structures. 
Moreover, if we regard $\{ Q_{log}^{\nu, i}( \mathscr N) \}_i$ as a cosimplicial pro-module on  $\{ \mathscr P_{X, log}^{\nu, i} \}_i$ (by restriction of scalars, via $\{ \alpha_{log}^{\nu,i} \}_i $), we see that it is the cosimplicial pro-module associated with the $\mathcal O_X-$pro-module with canonical stratification $\{ Q_{log}^{0, i}( \mathscr N) \}_i= \{ \mathscr P_{X,log}^{1,i} \}_i \otimes_{\mathcal O_X}  \mathscr N $. Indeed, for each $\nu \leq \mu$, $\{ Q_{log}^{\nu, i}( \mathscr N) \}_i$ is obtained from $\{ Q_{log}^{\mu, i}( \mathscr N) \}_i$ by base change with respect to any of the canonical morphisms $ \{ \mathscr P_{X, log}^{\nu, i} \}_i \longrightarrow \{ \mathscr P_{X, log}^{\mu, i} \}_i$. Now, for each  integer $k \geq 2$, we consider the differential operator $d^k$ of the log De Rham complex, 
$$
d^k \colon  \omega^k_X \longrightarrow  \omega^{k+1}_X
$$
Following \cite[p. 347]{Groth1}, for each $\nu \geq 0$,   $d^k$ induces a homomorphism of Artin-Rees pro-objects 
\begin{equation}
\label{o6}
Q_{log}^{\nu,^.}(d^k) \colon Q_{log}^{\nu,^.}( \omega^k_X ) \longrightarrow Q_{log}^{\nu,^.}(\omega^{k+1}_X)
\end{equation}
and we get the following cosimplicial complex   of  Artin-Rees pro-objects, 
\begin{equation}
\label{o7}
\{ Q_{log}^{*,i}(\mathcal O_X)  \}_i \longrightarrow \{ Q_{log}^{*,i}(\omega^1_X) \}_i  \longrightarrow 
\{ Q_{log}^{*,i}(\omega^2_X) \}_i  \longrightarrow ... \longrightarrow \{ Q_{log}^{*,i}(\omega^{k}_X)\}_i \longrightarrow ...
\end{equation}
The cosimplicial complex of Artin-Rees pro-objects $\{ Q_{log}^{*,i}(\omega^{^.}_X) \}_i $ 
is a resolution of $\omega^{^.}_X$ (\v{C}ech resolution). Indeed, we consider the double complex  of  $\mathcal O_X-$pro-modules 
\begin{equation}
\label{matrice0}
\omega^{^.}_X \stackrel{d_0}{\longrightarrow} \{ Q_{log}^{0,i}(\omega^{^.}_X) \}_i \stackrel{d_1- d_0}{\longrightarrow}  
\{ Q_{log}^{1,i}(\omega^{^.}_X) \}_i \stackrel{d_2- d_1 + d_0}{\longrightarrow}  ...
\longrightarrow \{ Q_{log}^{\nu,i}(\omega^{^.}_X) \}_i  \longrightarrow ...
\end{equation}
where the maps are obtained from  the cosimplicial maps \eqref{o2} (with respect to the cosimplicial index $\nu$), 
by ``forgetting one face" (\cite[p. 12]{Chamb}). 
Then, one can show that \eqref{matrice0} is locally  homotopic to zero, by using the degenerating 
maps of the cosimplicial complex  \eqref{o2} (\cite[\S V, Lemma 2.2.1]{Berth1}). 
Now, we apply to \eqref{matrice0} the additive functor $\{ \mathcal O_X/ \mathscr I^n \}_{n \in \mathbb N} \otimes_{\mathcal O_X} (-)$, 
in the category of pro-coherent $\mathcal O_X-$modules. Since it respects the local homotopies, we find that the complex
\begin{equation}
\label{matrice01}
\{ \mathcal O_X/ \mathscr I^n \}_{n} \otimes \omega^{^.}_X \stackrel{d_0}{\longrightarrow}
 \{ \mathcal O_X/ \mathscr I^n \}_{n} \otimes \{ Q_{log}^{0,i}(\omega^{^.}_X) \}_i \stackrel{d_1- d_0}{\longrightarrow} 
 \{ \mathcal O_X/ \mathscr I^n \}_{n } \otimes \{ Q_{log}^{1,i}(\omega^{^.}_X) \}_i \stackrel{d_2- d_1 + d_0}{\longrightarrow}  ... 
\end{equation}
is also locally homotopic to zero. \\ 
We give now a sort of ``Log Poincar\'e Lemma" in characteristic zero. 
\begin{theorem}
\label{teo0}
The complex of Artin-Rees $\mathcal O_X-$pro modules 
\begin{equation}
\label{2}
[\mathcal O_X \stackrel{d_0}{\longrightarrow} \{ L_{X}(\omega^{^.}_X) \}_i]=[\mathcal O_X \stackrel{d_0}{\longrightarrow}  \{ \mathscr P^{1,i}_{X,log} \}_i \longrightarrow  \{ \mathscr P^{1,i}_{X,log}   \}_i \otimes_{\mathcal O_X} \omega^{1}_X  \longrightarrow  \{ \mathscr P^{1,i}_{X,log}   \}_i \otimes_{\mathcal O_X} \omega^{2}_X \longrightarrow ...] 
\end{equation}
is locally homotopic to zero. 
\end{theorem}
\begin{Proof}
From  \eqref{equaz0}, it follows that the composition $\mathcal O_X \longrightarrow \{ L_{X}(\mathcal O_X)_i \}_i \longrightarrow \{  L_{X}(\omega^{1}_X)_i \}_i $ is zero, so \eqref{2} is in fact a  complex of $\mathcal O_X-$pro modules. It is easy to show that the complexes 
\begin{equation}
\label{eq22}
[\mathcal O_X \stackrel{d_0}{\longrightarrow}  \mathscr P^{1,i}_{X,log} \stackrel{L_X(d^{0}_{X})_i}{\longrightarrow}   \mathscr P^{1,i-1}_{X,log} \otimes_{\mathcal O_X} \omega^{1}_X  \stackrel{L_X(d^{1}_{X})_{i-1}}{\longrightarrow}  \mathscr P^{1,i-2}_{X,log} \otimes_{\mathcal O_X} \omega^{2}_X \stackrel{L_X(d^{2}_{X})_{i-2}}{\longrightarrow}  ...]
\end{equation}
are locally homotopic to zero, for each $i \in \mathbb N$. Indeed, from Proposition \ref{c01},  we can define the following maps on the local basis,
\begin{scriptsize}
$$
\begin{matrix}
\mathscr P^{1,i-p}_{X,log} \otimes_{\mathcal O_X} \omega^{p}_X & \stackrel{s_p}{\longrightarrow} &  \mathscr P^{1,i-p+1}_{X,log} \otimes_{\mathcal O_X} \omega^{p-1}_X \\
\xi_{1,i-p}^{\alpha_1} \cdots \xi_{r,i-p}^{\alpha_r} \otimes \xi_{i_1,1} \wedge \cdots \wedge \xi_{i_p,1} & \longmapsto & \frac{1}{k+p} \sum_{m=1}^p (-1)^{m+1} \xi_{1,i-p+1}^{\alpha_1} \cdots \xi_{r,i-p+1}^{\alpha_r} \xi_{i_m,i-p+1} \otimes \xi_{i_1,1} \wedge \cdots \wedge \hat{\xi}_{i_m, 1} \wedge \cdots \wedge \xi_{i_p,1} \\
\mathbb I \otimes \xi_{i_1,1} \wedge \cdots \wedge \xi_{i_p,1} & \longmapsto &   \frac{1}{p} \sum_{m=1}^p (-1)^{m+1}  \xi_{i_m,i-p+1} \otimes \xi_{i_1,1} \wedge \cdots \wedge \hat{\xi}_{i_m,1} \wedge \cdots \wedge \xi_{i_p,1}
\end{matrix}
$$
\end{scriptsize}
where $0 \neq \alpha_1+ \cdots + \alpha_r=k \leq i-p$, and extend them by linearity. It is easy to compute that $s_0 \circ d_0 =id$, and $L_X(d^{p-1}_{X})_{i-p+1} \circ s_p + s_{p+1} \circ  L_X(d^{p}_{X})_{i-p} = id$, for each $p \geq 0$ (see also \cite[Theorem 2.2]{Fiorot} for an alternative proof, in the classical case).
\end{Proof}\\
Now, since the additive functor $\{ \mathcal O_X/ \mathscr I^n \}_n \otimes_{\mathcal O_X} (-)$ respects the local homotopies, by Theorem \ref{teo0}, 
the following complex of  Artin-Rees $\mathcal O_X-$pro modules 
\begin{equation}
\label{eq222}
[\{ \mathcal O_X/ \mathscr I^n   \}_n \stackrel{d_0}{\longrightarrow}  \{   \mathcal O_X/ \mathscr I^n \}_n  \otimes  \{ \mathscr P^{1,i}_{X,log}  \}_i \longrightarrow  \{  \mathcal O_X/ \mathscr I^n \}_n \otimes  \{  \mathscr P^{1,i}_{X,log}   \}_i \otimes_{\mathcal O_X} \omega^{1}_X   \longrightarrow ... ] 
\end{equation}
is also locally homotopic to zero, in the category of pro-coherent $\mathcal O_X-$modules.  \\
We consider now the following double complex $(\star \star)$
$$
\begin{scriptsize} 
\begin{matrix}
  &  \{ \mathcal O_X/ \mathscr I^n \}_n & \longrightarrow \quad \quad \quad \quad  \{ \mathcal O_X/ \mathscr I^n \}_n   \otimes \omega^1_X  &  \longrightarrow \quad \quad \quad \quad \{ \mathcal O_X/ \mathscr I^n \}_n  \otimes \omega^2_X  &  ...\\
 &  ^{d_0}{\downarrow} & ^{d_0}{\downarrow} &  ^{d_0}{\downarrow} & \\
 \{ \mathcal O_X/ \mathscr I^n \}_n & \stackrel{d_0}{\longrightarrow}  \{  \mathcal O_X/ \mathscr I^n \}_n  \otimes  \{ \mathscr P^{1,i}_{X,log}   \}_i & \longrightarrow   
\{  \mathcal O_X/ \mathscr I^n \}_n  \otimes   \{ \mathscr P^{1,i}_{X,log}  \}_i \otimes  \omega^1_X  &  \longrightarrow  
\{  \mathcal O_X/ \mathscr I^n \}_n  \otimes   \{ \mathscr P^{1,i}_{X,log}  \}_i \otimes \omega^2_X   &  ...\\
^{d_1-d_0}{\downarrow}  &  ^{d_1-d_0}{\downarrow} &  ^{d_1-d_0}{\downarrow} &  ^{d_1-d_0}{\downarrow} & \\
\{  \mathcal O_X/ \mathscr I^n \}_n  \otimes  \{ \mathscr P^{1,i}_{X,log}   \}_i  &  \stackrel{d_0}{\longrightarrow}  \{  \mathcal O_X/ \mathscr I^n \}_n  \otimes  \{ \mathscr P^{2,i}_{X,log}   \}_i & \longrightarrow  \displaystyle 
\{   \mathcal O_X/ \mathscr I^n  \}_n  \otimes  \{  \mathscr P^{2,i}_{X,log}  \}_i \otimes \omega^1_X   &  \longrightarrow  
\{  \mathcal O_X/ \mathscr I^n \}_n   \otimes  \{ \mathscr P^{2,i}_{X,log}   \}_i \otimes \omega^2_X    &  ...\\
^{d_2-d_1+d_0}{\downarrow} &  ^{d_2-d_1+d_0}{\downarrow} &  ^{d_2-d_1+d_0}{\downarrow} &  ^{d_2-d_1+d_0}{\downarrow} & \\
... &  ... &  ... & ... & \\
\quad  \quad \quad \quad \downarrow &  \quad  \quad \quad \quad \downarrow &  \quad  \quad \quad \quad \downarrow &  \quad  \quad \quad \quad \downarrow & \\
\{  \mathcal O_X/ \mathscr I^n \}_n   \otimes \{  \mathscr P^{\nu-1,i}_{X,log}    \}_i &\stackrel{d_0}{\longrightarrow}   \{  \mathcal O_X/ \mathscr I^n \}_n   \otimes \{  \mathscr P^{\nu,i}_{X,log}    \}_i & \longrightarrow   
\{   \mathcal O_X/ \mathscr I^n \}_n  \otimes  \{   \mathscr P^{\nu,i}_{X,log}  \}_i \otimes  \omega^1_X   &  \longrightarrow 
\{   \mathcal O_X/ \mathscr I^n \}_n  \otimes  \{  \mathscr P^{\nu,i}_{X,log}   \}_i \otimes \omega^2_X   &  ...\\
\quad  \quad \quad \quad \downarrow &  \quad  \quad \quad \quad \downarrow &  \quad  \quad \quad \quad \downarrow &  \quad  \quad \quad \quad \downarrow & \\
... &  ... &  ... &  ... & \\
\end{matrix}
\end{scriptsize}
$$
Now, from \eqref{matrice01}, all the columns of $(\star \star)$, except the first, are locally homotopic to zero. Moreover, from \eqref{eq222}, 
the second row of $(\star \star)$ is also locally homotopic to zero.  The $(\nu+1)-$th row of this double complex ($\nu \geq 2$) 
is obtained from the second 
row by tensorizing (over $\mathcal O_X$) with the (log stratified) Artin-Rees pro-object $\{ \mathscr P^{\nu-1,i}_{X,log} \}_i$. Indeed, for 
each $\nu \geq 0$,  $\{ \mathscr P^{\nu,i}_{X,log} \}_i  \cong \{ \mathscr P^{\nu-1,i}_{X,log} \}_i \otimes \{ \mathscr P^{1,j}_{X,log} \}_j $. 
So, since the second row is locally homotopic to zero and the additive functor $\{ \mathscr P^{\nu-1,i}_{X,log} \}_i \otimes_{\mathcal O_X } (-)$ 
respects the local homotopies, 
we see that each row of $(\star \star)$, except the first,  is  also locally homotopic to zero. \\ 
Therefore, we can conclude that  the double complex $\{  \mathcal O_X/ \mathscr I^n \}_n \otimes \{  Q^{*,i}_{log} (\omega^{^.}_X) \}_i $ 
is a resolution of both the first column $\{ \mathcal O_X/ \mathscr I^n\}_n  \otimes \{ \mathscr P^{*,i}_{X,log}  \}_i$,  
and the first row $ \{ \mathcal O_X/ \mathscr I^n \}_n   \otimes \omega^{^.}_X$ of $(\star \star)$. 
Then, since all pro-systems satisfy the Mittag-Leffler condition, we get the two following isomorphisms in cohomology,  
\begin{equation}
\label{eq224}
\mathbb H^{^.}(Y,  \displaystyle \lim_{\stackrel{\longleftarrow}{n}} \lim_{\stackrel{\longleftarrow}{i}} \mathcal O_X/ \mathscr I^n \otimes \mathscr P^{*,i}_{X,log} ) \cong \mathbb H^{^.}(Y, \displaystyle \lim_{\stackrel{\longleftarrow}{n}} \lim_{\stackrel{\longleftarrow}{i}} \mathcal O_X/ \mathscr I^n \otimes  Q^{*,i}_{log} (\omega^{^.}_X) )
\end{equation}
\begin{equation}
\label{equazione1}
H^{^.}_{DR,log}((Y,N)/\mathbb C):=\mathbb H^{^.}(Y, \displaystyle  \lim_{\stackrel{\longleftarrow}{n}} \mathcal O_X/ \mathscr I^n    \otimes \omega^{^.}_X ) \cong \mathbb H^{^.}(Y, \displaystyle \lim_{\stackrel{\longleftarrow}{n}} \lim_{\stackrel{\longleftarrow}{i}} \mathcal O_X/ \mathscr I^n    \otimes Q_{log}^{*,i} (\omega^{^.}_X))
\end{equation}
\begin{remark}
\label{oss1}
{\rm Since the Artin-Rees $\mathcal O_X-$pro module $\{ \mathcal O_X/ \mathscr I^n   \}_n$ is endowed with a log stratification (see Remark \ref{re00}), we have isomorphisms, for any $\nu, k \geq 0$,  
$$
\{ \mathcal O_X/ \mathscr I^n   \}_n \otimes \{   \mathscr P^{\nu,i}_{X,log} \}_i \otimes \omega^{^.}_X  \cong \{   \mathscr P^{\nu,i}_{X,log} \}_i  \otimes \{ \mathcal O_X/ \mathscr I^n   \}_n \otimes \omega^{^.}_X 
$$
and so there is an identification 
$$
\{  \mathcal O_X/ \mathscr I^n \}_n \otimes \{  Q^{*,i}_{log} (\omega^{^.}_X) \}_i \cong \{  Q^{*,i}_{log} (\{  \mathcal O_X/ \mathscr I^n \}_n \otimes \omega^{^.}_X) \}_i 
$$}
\end{remark}

In order to calculate $H^{^.}(Y^{log}_{inf}, \mathcal O_{Y^{log}_{inf}})$, from the exact closed 
imersion $Y \hookrightarrow X$ one can define the sheaf $\displaystyle \tilde{Y}:=   \lim_{\stackrel{\longrightarrow}{n}} \tilde{Y}_n$ on $InfLog(Y/S)$ 
($Y_n$ being the $n-$th log infinitesimal neighborhood of $Y$ in $X$), which covers the final object  of  $Y^{log}_{inf}$ (\cite[\S 5.2]{Groth1}). 
If $(X^{\nu} ) \hat{|} Y$ is the log formal completion of $X^{\nu}$ along its closed 
log subscheme $Y \hookrightarrow X \hookrightarrow X^{\nu}$, we can introduce the (\v{C}ech) cosimplicial 
sheaf on $Y$ equal to $\mathcal O_{(X^{*}) \hat{|} Y}$ 
(this is denoted by $\mathscr F^{*}$ in \cite[p. 338]{Groth1}, when $F=\mathcal O_{Y^{log}_{inf}}$). Then, by means of $\tilde{Y}$ 
and this \v{C}ech cosimplicial sheaf,  we 
have that 
$\mathbb H^{^.}(Y, \mathcal O_{(X^{*}) \hat{|} Y}) \cong H^{^.}(Y^{log}_{inf}, \mathcal O_{Y^{log}_{inf}})$ (\cite[\S 5.1 and p. 339]{Groth1}). 
Finally, since 
$\displaystyle \lim_{\stackrel{\longleftarrow}{n}} \lim_{\stackrel{\longleftarrow}{i}} \mathcal O_X/ \mathscr I^n \otimes \mathscr P^{*,i}_{X,log}$, 
as  cosimplicial sheaf on $Y$, is equal to $\mathcal O_{(X^{*}) \hat{|} Y}$, \begin{equation}
\label{eq225}
\mathbb H^{^.}(Y, \displaystyle \lim_{\stackrel{\longleftarrow}{n}} \lim_{\stackrel{\longleftarrow}{i}} \mathcal O_X/ \mathscr I^n \otimes \mathscr
 P^{*,i}_{X,log} ) \cong H^{^.}(Y^{log}_{inf}, \mathcal O_{Y^{log}_{inf}})
\end{equation}
From  \eqref{eq224}, \eqref{equazione1} and \eqref{eq225}, we conclude that there exists an isomorphism between the log De Rham cohomology of $Y$ and its log Infinitesimal cohomology, namely
\begin{equation}
\label{eq226}
H^{^.}_{DR,log}((Y,N)/\mathbb C) \cong H^{^.}(Y^{log}_{inf}, \mathcal O_{Y^{log}_{inf}})
\end{equation}
\begin{remark}
\label{eq2206}
Forgetting log structure, one could use  analogous techniques in the classical setting (which are nothing but a miscellanea of those ones indicated in \cite{Groth1}) and obtain the isomorphism \eqref{eq226} for any scheme $Y$ over $S$ (without log). This fact, together with the result proved by 
Hartshorne in \cite[Chapter IV, Theorem I.I]{Hartshorne}, gives the isomorphisms between \eqref{01}, \eqref{02} and \eqref{03} in the Introduction. 
\end{remark}

\section{The Complex $\omega^{^.,log}_{X \hat{|}Y}$}
\label{o}

With $X$ and $Y$ as in \S 1, let $Y^{an}$, $X^{an}$ be the associated fs log analytic spaces, and let
 $i^{an} \colon Y^{an } \hookrightarrow X^{an} $ be the corresponding analytic exact closed immersion. 
When the context obviates any confusion, we will omit the superscript $(-)^{an}$ in denoting the associated analytic spaces.  \\
We consider the closed analytic subspaces $Y_k$  of $X$, defined by the ideals $\mathscr I^k$, with $k \in \mathbb N$. 
 On each such $Y_k$ we consider the log structure induced by $M_X$, i.e., if $i_k \colon Y_k \hookrightarrow X$  is  the closed immersion, then we take $M_{Y_k} = i_k^* M_X$. We have a sequence of exact closed immersions, which we denote by $\varphi_k$,
$$
Y=Y_1 \stackrel{\varphi_1}{\hookrightarrow} Y_2  \stackrel{\varphi_2}{\hookrightarrow} Y_3 \stackrel{\varphi_3}{\hookrightarrow}... \stackrel{\varphi_k}{\hookrightarrow} Y_{k+1} \stackrel{\varphi_{k+1}}{\hookrightarrow} ... \hookrightarrow X
$$
Therefore, we have 
a projective system of rings $\{ \mathcal O_{Y_k} \cong i_k^{-1}(\mathcal O_{X}/ \mathscr I^k); \, \varphi_k \colon \mathcal O_{Y_{k+1}} \longrightarrow \mathcal O_{Y_{k}} \}_{k \geq 1}$, where the transition  maps  $\varphi_k $ are surjective. Moreover, the diagram
\begin{equation}
\label{re0}
\begin{matrix}
\varphi_k^{-1}(M_{Y_{k+1}}) & \longrightarrow  &  M_{Y_{k}} \\
^{\alpha_{k+1}}\downarrow & & ^{\alpha_k}{\downarrow}  \\
\varphi_k^{-1}( \mathcal O_{Y_{k+1}} )& \longrightarrow  & \mathcal O_{Y_k}
\end{matrix}
\end{equation}
is commutative, for each $k \geq 1$. Since $M_X$ is a fine log structure on $X$,  each $Y_k$ is endowed with a fine log structure $i_k^* M_X$, which
we donote more  by $M_k$.\\
The closed immersion $i_k \colon Y_k \hookrightarrow X$ is exact, for each $k$, so (\cite[(1.4.1)]{Kato1})
\begin{equation}
\label{re1}
M_k/\mathcal O_{Y_k}^* =  i_k^* M_X/\mathcal O_{Y_k}^* \cong i_k^{-1} (M_X/\mathcal O_X^*) = (M_X/\mathcal O_X^*)_{|Y_k}
\end{equation}
\begin{remark}
\label{re2}
Since the underlying topological space of each $Y_k$ is equal to $Y$, it follows from \eqref{re1} that 
$M_k/\mathcal O_{Y_k}^*  \cong (M_X/\mathcal O_X^*)_{|Y}$, for each $k \geq 1$. 
Therefore, if we consider the associated sheaf of groups $(M_k/\mathcal O_{Y_k}^* )^{gp} \cong M^{gp}_k/\mathcal O_{Y_k}^* $, 
we have that, for each $k \geq 1$,
\begin{equation}
\label{re3}
M^{gp}_k/\mathcal O_{Y_k}^*  \cong (M^{gp}_X/\mathcal O_X^*)_{|Y}
\end{equation}
\end{remark}

Let $Y^{log}_k$ (resp. $X^{log}$) be the Kato-Nakayama topological space associated to $Y_k$ (resp. to $X$), and let $\tau_k \colon Y^{log}_k \longrightarrow Y_k$ (resp. $\tau_X \colon X^{log} \longrightarrow X$) be the corresponding  surjective, continuous and proper map of topological spaces (\S \ref{Weight0}). 
We now consider  the ``formal" analytic space $X \hat{|} Y$,  
which is $Y^{an}$ as topological space, and whose structural sheaf is 
$$
\mathcal O_{X \hat{|} Y} =: \lim_{\stackrel{\longleftarrow}{k}} \mathcal O_{Y_k} \cong  \lim_{\stackrel{\longleftarrow}{k}} i_k^{-1}(\mathcal O_{X}/ \mathscr I^k)
$$ 
Now, since the closed immersion $i \colon Y \hookrightarrow X$ is exact, the formal completion of the fs log analytic space $X$ along the closed log subspace $Y$ is equal to the classical completion $X \hat{|} Y$, endowed with the log structure induced by $M_X$ (Definition \ref{def1}). So, if $i_{X \hat{|} Y} \colon X \hat{|} Y \hookrightarrow X$, then the log structure on $X \hat{|} Y$ is $i_{X \hat{|} Y}^* M_X$.  We denote it by $M_{X \hat{|}Y}$. 
We now  define a ringed topological space $((X \hat{|}Y)^{log}, \mathcal O^{log}_{X \hat{|}Y})$, associated to the formal fine log analytic space $X \hat{|}Y$. 
\begin{definition}
\label{re4}
With the previous notation, we define $(X \hat{|}Y)^{log}$ to be the topological space $Y^{log}$, endowed with the following  sheaf of rings
\begin{equation}
\label{re5}
\mathcal O^{log}_{X \hat{|}Y}=: \tau_Y^{-1}(\mathcal O_{X \hat{|}Y}) \otimes_{\tau_X^{-1} (\mathcal O_X)} (\mathcal O^{log}_{X})
\end{equation}
\end{definition}

\begin{Lemma}\cite[Lemma (3.3)]{Nakayama}
\label{re6}
With the previous notation, let $x \in Y$, $y \in Y^{log}$ be such that $\tau_Y(y)=x$. Let $\mathscr L_X$ be the sheaf of logarithms 
of local sections of $\tau_X^{-1}(M_X^{gp})$ (\S \ref{Weight0}). Let $\{ t_1,...,t_n\}$ be a family of elements of the stalk $\mathscr L_{X,y}$,
 whose image under the map ${\rm exp}_y \colon \mathscr L_{X,y} \longrightarrow \tau_X^{-1}(M^{gp}_X/\mathcal O_X^*)_y$ is a $\mathbb Z-$basis 
of $M^{gp}_{X,x}/\mathcal O_{X,x}^*$. Then, $\mathcal O^{log}_{(X \hat{|}Y),y}$ is isomorphic, as $\mathcal O_{(X \hat{|}Y),x}-$algebra, 
to the polynomial ring $\mathcal O_{(X \hat{|}Y),x}[T_1,...,T_n]$, via the correspondence 
\begin{equation}
\label{re06}
\begin{matrix}
\mathcal O_{(X \hat{|}Y),x}[T_1,...,T_n]  & \longrightarrow  & \mathcal O^{log}_{(X \hat{|}Y),y}\\
 & & \\
T_i & \longmapsto & t_i
\end{matrix}
\end{equation}
for $i=1,...,n$. 
\end{Lemma}
\begin{Proof}
By \cite[Lemma (3.3)]{Nakayama}, applied to $X^{log}$, the isomorphism 
$\mathcal O^{log}_{X,y} \cong \tau_X^{-1}(\mathcal O_X)_y [T_1,...,T_n]$ implies that 
$\mathcal O^{log}_{(X \hat{|}Y),y} 
\cong \tau_Y^{-1} ( \mathcal O_{X \hat{|}Y})_y \otimes_{\tau_X^{-1} (\mathcal O_X)_y} \tau_X^{-1}(\mathcal O_X)_y[T_1,...,T_n]  
\cong 
\tau_Y^{-1}(\mathcal O_{X \hat{|}Y})_y[T_1,...,T_n] \cong$  \\ $\mathcal O_{(X \hat{|}Y),x}[T_1,...,T_n]$. \end{Proof}
\begin{Lemma}
\label{re8}\cite[Lemma (3.4)]{Nakayama}  
Let $r \in \mathbb Z$. We define a filtration $\hat{{\rm fil}}_r(\mathcal O^{log}_{X\hat{|}Y})$ on $\mathcal O^{log}_{X\hat{|}Y}$ by
\begin{equation}
\label{re08}
\hat{{\rm fil}}_r(\mathcal O^{log}_{X\hat{|}Y})=: \tau^{-1}(\mathcal O_{X \hat{|}Y}) \otimes_{\tau^{-1}(\mathcal O_X)} {\rm fil}_r(\mathcal O^{log}_{X})
\end{equation}
(where ${\rm fil}_r(\mathcal O^{log}_{X})$ is defined by Kato-Nakayama as 
${\rm Im} \{ \tau^{-1}(\mathcal O_X) \otimes_{\mathbb Z} (\bigoplus_{j=1}^r 
{\rm Sym}^j_{\mathbb Z} \mathscr L_X) \longrightarrow \mathcal O^{log}_X \}$). Then, the canonical map 
$$
\tau^{-1}(M^{gp}_X/\mathcal O^*_X) \cong \mathscr L_X/ \tau^{-1}(\mathcal O_X) \subseteq {\rm fil}_1(\mathcal O^{log}_{X})/ {\rm fil}_0(\mathcal O^{log}_{X})
$$
 induces the following  isomorphism
\begin{equation}
\label{re008}
\tau^{-1}(\mathcal O_{X \hat{|}Y}) \otimes_{\mathbb Z} \tau^{-1}({\rm Sym}^r_{\mathbb Z} (M^{gp}_X/\mathcal O^*_X))
\stackrel{\cong}{\longrightarrow}  \hat{{\rm fil}}_r(\mathcal O^{log}_{X\hat{|}Y})/\hat{{\rm fil}}_{r-1}(\mathcal O^{log}_{X\hat{|}Y})
\end{equation}
\end{Lemma}
\begin{Proof}
By \cite[Lemma (3.4)]{Nakayama}, for any $r \geq 0$, we have an isomorphism
\begin{equation}
\label{re0008}
\tau^{-1}(\mathcal O_X) \otimes_{\mathbb Z} \tau^{-1}({\rm Sym}^r_{\mathbb Z} (M^{gp}_X/\mathcal O^*_X)) 
\stackrel{\cong}{\longrightarrow}  {\rm fil}_r(\mathcal O^{log}_{X})/{\rm fil}_{r-1}(\mathcal O^{log}_{X})
\end{equation}
So, 
$$
\tau^{-1}(\mathcal O_{X \hat{|}Y}) \otimes_{\mathbb Z} \tau^{-1}({\rm Sym}^r_{\mathbb Z} (M^{gp}_X/\mathcal O^*_X)) \cong 
$$
$$
\tau^{-1}(\mathcal O_{X \hat{|}Y}) 
\otimes_{\tau^{-1}(\mathcal O_X)} 
\left( \tau^{-1}(\mathcal O_X) \otimes_{\mathbb Z} \tau^{-1}({\rm Sym}^r_{\mathbb Z} (M^{gp}_X/\mathcal O^*_X)) \right)  
$$
and, by \eqref{re0008}, this is isomorphic to 
\begin{equation}
\label{re00008}
\tau^{-1}(\mathcal O_{X \hat{|}Y}) \otimes_{\tau^{-1}(\mathcal O_X)} {\rm fil}_r(\mathcal O^{log}_{X})/{\rm fil}_{r-1}(\mathcal O^{log}_{X})
\end{equation}
Now, since the functor $ \tau^{-1}(\mathcal O_{X \hat{|}Y}) \otimes_{\tau^{-1}(\mathcal O_X)} ( -) $ is right exact, it follows that \eqref{re00008} is
isomorphic  to $\hat{{\rm fil}}_r(\mathcal O^{log}_{X\hat{|}Y})/\hat{{\rm fil}}_{r-1}(\mathcal O^{log}_{X\hat{|}Y})$. 
\end{Proof}

Using the  Kato-Nakayama complex $\omega^{^.,log}_X$ associated
 to the fs log smooth log analytic space $X$ over $S$ (\cite[(3.5)]{Nakayama}), we can now give the following
\begin{definition}
\label{re9}
In the previous notation, for any $q \in \mathbb N$, $0 \leq q \leq {\rm rk}_{\mathbb Z} \, \omega^1_X$, we define the following sheaf on $Y^{log}$
\begin{equation}
\label{re09}
\omega^{q,log}_{X \hat{|}Y}=: \mathcal O^{log}_{X \hat{|}Y} \otimes_{\tau^{-1}(\mathcal O_X)} \tau^{-1}(\omega^q_X)
\end{equation}
\end{definition}
Since $X$ is log smooth over $S$, it follows that $\omega^q_X$ is a locally free $\mathcal O_X-$module
 of finite type, and so $\omega^{q,log}_{X \hat{|}Y}$ is a locally free $\mathcal O^{log}_{X \hat{|} Y}-$module of finite type. 
Moreover, since $\omega^{q,log}_{X}=: \mathcal O^{log}_{X} \otimes_{\tau^{-1}(\mathcal O_X)} \tau^{-1}(\omega^q_X)$ (\cite[(3.5)]{Nakayama}),
 using the definition of $\mathcal O^{log}_{X \hat{|}Y} $ \eqref{re5}, we can also write 
\begin{equation}
\label{re10}
\omega^{q,log}_{X \hat{|}Y}= \tau^{-1}(\mathcal O_{X \hat{|}Y}) \otimes_{\tau^{-1}(\mathcal O_X)} \omega^{q,log}_X
\end{equation}
Now, since the differential $d^q \colon \omega^{q,log}_X \longrightarrow \omega^{q+1,log}_X$ is induced by that of $\omega^{^.}_X$, 
we have that $d^1 (\mathscr I^r) \subseteq \mathscr I^{r-1}\omega^1_X$, for any $r \geq 0$. Thus, $d^{^.} $ can be extended to get a differential
\begin{equation}
\label{re11}
\hat{d}^q \colon \omega^{q,log}_{X \hat{|}Y} \longrightarrow \omega^{q+1,log}_{X \hat{|}Y}
\end{equation}
Therefore, we obtain a complex $\omega^{^., log}_{X \hat{|}Y}$, whose differentials $\hat{d}^{^.}$ satisfy $\hat{d}^1(x)={\rm dlog}({\rm exp} (x))$,
 for each element $x \in \mathscr L_X$, and $\hat{d}^1(y) \in \mathscr I^{r-1}\omega^1_X$, for each element $y \in \mathscr I^r$, $r \geq 0$.

\section{Formal Poincar\'e Residue Map}
\label{formalres}

In this section, we want to give a ``formal version" of the Poincar\'e Residue map given by Deligne (\cite[(3.1.5.2)]{Deligne2}). 
We consider an fs log scheme $Y$, with log structure $M_Y$, and an exact closed immersion $i  \colon Y \hookrightarrow X$, 
where $X$ is an fs log smooth log scheme, 
with log structure $M_X$. We also suppose that the underlying scheme of $X$ is smooth over $S$, and 
its log structure $M_X$ is given by a normal crossing divisor 
$D \hookrightarrow X$, i.e. $M_X= j_* \mathcal O_U^* \cap \mathcal O_X \hookrightarrow \mathcal O_X$, 
where $j \colon U=X-D \hookrightarrow X$ is the open immersion. Let $X^{an}$, $Y^{an}$ be the log analytic spaces associated 
to $X$ and $Y$, which we will simply denote by $X$, $Y$,
 when no  confusion can arise. \\
We take the log De Rham complex $\omega^{^.}_X= \Omega^{^.}_X({\rm log} \, M_X)= \Omega^{^.}_X({\rm log} \, D)$. 
Its completion  $\omega^{^.}_{X \hat{|} Y}$,  along the closed subscheme $Y$ of $X$, satisfies   
$$
\omega^{i}_{X \hat{|}Y} \cong \omega^i_X \otimes_{\mathcal O_X} \mathcal O_{X \hat{|}Y}
$$
for each $i$, $0 \leq i \leq n={\rm dim} \, X$, because the $\mathcal O_X-$modules $\omega^i_X$ are locally free, and so coherent. \\
We denote by $\hat{d}^i$ the differential $\omega^i_{X \hat{|}Y} \longrightarrow \omega^{i+1}_{X \hat{|}Y}$ of the complex
 $\omega^{^.}_X \otimes_{\mathcal O_X} \mathcal O_{X \hat{|}Y}$. 
We consider the  weight filtration $W.$ on  $\omega^{^.}_X$ (\cite[(3.1.5.1)]{Deligne2}): since each term 
$W_k  (\omega^i_X)= \Omega^{i-k}_X \wedge \omega^k_X$ is a locally free $\mathcal O_X-$module, 
this filtration induces an increasing filtration $\hat{W}_.$ on $\omega^{^.}_{X \hat{|}Y}$, defined by 
$$
\hat{W}_k(\omega^{i}_{X \hat{|}Y})=: W_k (\omega^{i}_X) \otimes_{\mathcal O_X} \mathcal O_{X \hat{|}Y}
$$
for each $0\leq k \leq i \leq n$. \\
Since $\hat{W}_k (\omega^{i}_{X \hat{|}Y}) = {\rm Im} \{ \omega^k_X \otimes_{\mathcal O_X} \Omega^{i-k}_X \otimes_{\mathcal O_X} 
\mathcal O_{X \hat{|}Y} \stackrel{\wedge \otimes id}{\longrightarrow} \omega^i_X \otimes_{\mathcal O_X} \mathcal O_{X \hat{|}Y} \}$, we can write the 
term $\hat{W}_k$ of the filtration as  $\omega^k_X \wedge \Omega^{i-k}_{X\hat{|}Y}$, where $\Omega^{^.}_{X \hat{|}Y}$ is the completion 
of the classical 
De Rham complex $\Omega^{^.}_X$ along $Y$. Moreover, we note that $\hat{W}_k (\omega^{i}_{X \hat{|}Y})$ is a locally free
 $\mathcal O_{X \hat{|} Y}-$submodule 
of $\omega^{i}_{X \hat{|}Y}$, for each $i$.\\
We suppose now that, locally at a point $x \in  Y \hookrightarrow X$, the normal crossing divisor $D$ is the union of smooth irreducible 
components $D= D_1 \cup ... \cup D_r$, where each component $D_i$ is locally defined by the equation $z_i=0$ 
(for a local coordinate system $\{ z_1,...,z_n \}$ of $X$ at $x$). Let $S^k$ be the set of strictly increasing sequences of 
indices $\sigma= (\sigma_1,..., \sigma_k)$, where $\sigma_i \in \{ 1,...,r \}$, 
and let $D_{\sigma} = D_{\sigma_1} \cap ... \cap D_{\sigma_k}$. Let $D_k=\bigcup_{\sigma \in S^k} D_{\sigma}$ and $D^k$ be the 
disjoint union $\coprod_{\sigma \in S^k} D_{\sigma}$. Moreover, let $\pi^k \colon D^k \longrightarrow X$ be the canonical map. Then, locally at $x$,
 the $\mathcal O_{X \hat{|}Y}$-submodule $\hat{W}_k (\omega^{i}_{X \hat{|}Y})$ of $\omega^{i}_{X \hat{|}Y}$ can be written as
$$
\hat{W}_k (\omega^{i}_{X \hat{|}Y})= \sum_{\sigma \in  S^k} \Omega^{i-k}_{X \hat{|}Y} \wedge {\rm dlog} z_{\sigma_1} \wedge ... \wedge {\rm dlog} z_{\sigma_k}
$$
for each $0\leq i \leq n$. Therefore,  the elements of $\hat{W}_k (\omega^{i}_{X \hat{|}Y})$ are locally linear combinations of terms  
$\eta \wedge {\rm dlog} z_{\sigma_1} \wedge ... \wedge {\rm dlog} z_{\sigma_k}$, with $\eta \in \Omega^{i-k}_{X \hat{|}Y}$.\\
Let  $Y_k= D_k \cap Y$, and $Y^k=\coprod_{\sigma \in S^k} (Y_{\sigma})$, with $Y_{\sigma}= D_{\sigma} \cap Y$. 
We have the following cartesian diagram
\begin{equation}
\label{res0}
\begin{matrix}
Y^k & \hookrightarrow  & D^k  \\
^{\pi^k_Y}\downarrow & & ^{\pi^k}{\downarrow}  \\
Y & \stackrel{i}{\hookrightarrow}  & X 
\end{matrix}
\end{equation}
Since each intersection $D_{\sigma}$ is smooth over $S$, we can take the sheaf of classical differential $i- $forms $\Omega^i_{D^k}$ over $D^k$; then, $\pi^k_*(\Omega^i_{D^k}) \cong \bigoplus_{\sigma \in S^k} (i_{\sigma*} \Omega^i_{D_{\sigma}})$, where $i_{\sigma} \colon D_{\sigma} \hookrightarrow X$.  \\
So, $(\pi^k_*(\Omega^i_{D^k}))_{\hat{|}Y} \cong (\bigoplus_{\sigma \in S^k} (i_{\sigma*} \Omega^i_{D_{\sigma}}))_{\hat{|}Y} \cong 
\bigoplus_{\sigma \in S^k} (i_{\sigma*} \Omega^i_{D_{\sigma}})_{\hat{|}Y}$. \\
From the cartesian diagram 
\begin{equation}
\label{res1}
\begin{matrix}
Y_{\sigma} & \hookrightarrow  & D_{\sigma} \\
\downarrow & & ^{i_{\sigma}}{\downarrow}  \\
Y & \stackrel{i}{\hookrightarrow}  & X 
\end{matrix}
\end{equation}
we deduce the map $\hat{i}_{\sigma} \colon D_{\sigma} \hat{|} Y_{\sigma} \hookrightarrow X \hat{|} Y$. From \cite[(Corollaire (10.14.7))]{Groth2}, it follows that 
$$
(i_{\sigma*} \Omega^i_{D_{\sigma}})_{\hat{|}Y} \cong \hat{i}_{\sigma*}(\Omega^i_{D_{\sigma} \hat{|} Y_{\sigma}})
$$
and then 
\begin{equation}
\label{res2}
(\pi^k_*\Omega^i_{D^k})_{\hat{|}Y} \cong \bigoplus_{\sigma \in S^k} \hat{i}_{\sigma*}(\Omega^i_{D_{\sigma} \hat{|} Y_{\sigma}})
\end{equation}

\subsection{The Formal Poincar\'e Residue}
\label{res3}

In \cite[(3.1.5.2)]{Deligne2}, Deligne defines a map of complexes
\begin{equation}
\label{res03}
{\rm Res^{^.}} \colon {\rm Gr}^W_k(\Omega^{^.}_X({\rm log} \, D)) \longrightarrow \pi^k_* \Omega^{^.}_{D^k}(\varepsilon^k)[-k]
\end{equation}
for each $k \leq n$, called the Poincar\'e Residue map, where $\varepsilon^k$ is defined as in \cite[(3.1.4)]{Deligne2}, and represents the orientations of the intersections $D_{\sigma}$ of $k$ components of $D$. Given a local section $\eta \wedge {\rm dlog}z_{\sigma_1} \wedge... \wedge  {\rm dlog}z_{\sigma_k}  \in {\rm Gr}^W_k(\Omega^{p}_X({\rm log} \, D))$, with $\eta \in \Omega^{p-k}_X$, the map ${\rm Res}$ sends it to $\eta_{| D_{\sigma}} \otimes ({\rm orientation} \, \sigma_1... \sigma_k)$. Deligne proved that ${\rm Res}$ is an isomorphism of complexes (\cite{Deligne2}). Moreover, from \cite[(3.1.8.2)]{Deligne2},  the following sequence of isomorphisms 
$$
\mathbb R^k j_* \mathbb C \cong \mathscr H^k(j_* \Omega^{^.}_U) \cong \mathscr H^k(\Omega^{^.}_X({\rm log} \, D)) \cong \varepsilon^k_X
$$
implies that there exists an identification 
\begin{equation}
\label{res10}
\varepsilon^k_X \cong \mathbb C \otimes_{\mathbb Z} \bigwedge^k M^{gp}_X/\mathcal O_X^*
\end{equation}
($M_X$ is the log structure on $X$ associated to the normal crossing divisor $D$, and $\varepsilon^k_X$ is the direct 
image of $\varepsilon^k$ via the map $\pi^k \colon D^k \longrightarrow X$ \cite[(3.1.4.1)]{Deligne2}). 
Using diagram \eqref{res0}, and \eqref{res2}, we can extend the Deligne Poincar\'e Residue map to the formal case.  
\begin{definition}
\label{res4}
In the previous notation, we define the map
\begin{equation}
\label{res04}
\hat{{\rm Res}}^{p} \colon {\rm Gr}^W_k (\omega^p_{X \hat{|}Y}) \cong {\rm Gr}^W_k (\omega^p_X \otimes_{\mathcal O_X} 
\mathcal O_{X \hat{|} Y}) \longrightarrow \hat{\pi}^k_* \Omega^{p-k}_{D^k \hat{|} Y^k} (\varepsilon^k)
\end{equation}
as
$$
\eta \wedge {\rm dlog} z_{\sigma_1} \wedge ... \wedge {\rm dlog} z_{\sigma_k} \longmapsto \eta_{| (D_{\sigma} \hat{|} Y_{\sigma})}
 \otimes ({\rm orientation} \, \sigma_1... \sigma_k)
$$
where $\eta \in \Omega^{p-k}_{X \hat{|} Y}$. 
\end{definition}
We consider  the completion $\hat{{\rm Res}}^p$ of the Deligne Poincar\'e Residue map \eqref{res03}, in degree $p$, along the closed subscheme $Y$. 
We will prove that the maps $\hat{{\rm Res}}^p$ induce the following $\mathcal O_X-$ linear isomorphism of complexes
\begin{equation}
\label{res5}
\hat{{\rm Res}}^{^.} \colon {\rm Gr}^W_k (\omega^{^.}_{X \hat{|}Y})
 \stackrel{\cong}{\longrightarrow} \hat{\pi}^k_* \Omega^{^.}_{D^k \hat{|} Y^k} (\varepsilon^k)[-k]
\end{equation} 
for each $k \leq n$. \\
To this end, we briefly recall the classical construction of the Deligne Poincar\'e  Residue map. 
So, given $\sigma \in S^k$, we consider the application 
\begin{equation}
\label{rho}
\rho_{\sigma} \colon \Omega^{p-k}_X \longrightarrow {\rm Gr}^W_k(\omega^{p}_X)
\end{equation}
which is locally defined by 
\begin{equation}
\label{res05}
\rho_{\sigma}(\eta)=: \eta \wedge {\rm dlog} z_{\sigma_1} \wedge ... \wedge {\rm dlog} z_{\sigma_k}
\end{equation}
This map does not depend on the choice of the local coordinates $z_i$ (\cite[3.6.6]{Elzein}). Moreover we have that 
$$
\rho_{\sigma}(z_{\sigma_i} \cdot \beta)=0 \quad {\rm and } \quad \rho_{\sigma}(dz_{\sigma_i} \wedge \gamma)=0
$$
for all sections $\beta \in \Omega^{p-k}_X$, and $\gamma \in \Omega^{p-1-k}_X$.  Therefore $\rho_{\sigma}$ factorizes into 
\begin{equation}
\label{res6}
\begin{matrix}
\Omega^{p-k}_X & \longrightarrow & i_{\sigma*} \Omega^{p-k}_{D_{\sigma}} \otimes ({\rm orientation} \, \sigma_1... \sigma_k) \\
^{\rho_{\sigma}}{\downarrow} & {\swarrow}_{\overline{\rho}_{\sigma}} &  \\
 {\rm Gr}^W_k(\omega^{p}_X) & & 
\end{matrix}
\end{equation}
Thus, all these maps being locally compatible with the differentials, the maps $\overline{\rho}_{\sigma}$ define a morphism of complexes 
\begin{equation}
\label{res7}
\overline{\rho}^{^.} \colon \pi^k_* \Omega^{^.}_{D^k} (\varepsilon^k)[-k] \longrightarrow {\rm Gr}^W_k(\omega^{^.}_X)
\end{equation}
This morphism is locally defined by \eqref{res05}, and it is a global morphism on $X$: it is an isomorphism of complexes. 
Its inverse isomorphism $ {\rm Gr}^W_k(\omega^{^.}_X) \longrightarrow \pi^k_* \Omega^{^.}_{D^k}  (\varepsilon^k)[-k] $ 
is the Deligne Poincar\'e Residue map ${\rm Res}^{^.}$ (\cite[(3.6.7.1)]{Elzein}).  In view of this construction,  
we can see that ${\rm Res}^{^.}$ is an $\mathcal O_X-$linear morphism of complexes: so we deduce that the maps $\hat{{\rm Res}}^{p}$ 
in \eqref{res04}, are compatible with the differentials induced from $\omega^{^.}_{X \hat{|}Y}$  and $\Omega^{^.}_{D^k \hat{|} Y^k} [-k]$, 
because $\hat{{\rm Res}}^{p}$ comes from the $\mathcal O_X-$linear map ${\rm Res}^p$ by completion along $Y$. Indeed, we note that 
$$
 {\rm Gr}^W_k (\omega^{^.}_{X \hat{|}Y})  \cong  {\rm Gr}^W_k (\omega^{^.}_{X }) \otimes_{\mathcal O_X} \mathcal O_{X \hat{|}Y}
$$
and, from \eqref{res2}, we have that 
$$
(\pi^k_* \Omega^{^.}_{D^k } (\varepsilon^k)[-k])_{\hat{|}Y} \cong  \bigoplus_{\sigma \in S^k} \hat{i}_{\sigma*}(\Omega^{^.}_{D_{\sigma} \hat{|} Y_{\sigma}}) (\varepsilon^k)[-k] \cong \hat{\pi}^k_* \Omega^{^.}_{D^k \hat{|} Y^k} (\varepsilon^k)[-k]
$$
So we conclude that the morphism of complexes $\hat{{\rm Res}}^{^.}$ \eqref{res5} is an isomorphism, for each $k\leq n$.
\begin{remark}
\label{res8}
We can also construct the morphism $\hat{{\rm Res}}^{^.} $ using a formal version of the classical construction  of ${\rm Res}^{^.} $, described in \eqref{rho}, \eqref{res05}, \eqref{res6}, \eqref{res7}. Indeed, we can define the map 
\begin{equation}
\label{res08}
\rho_{\sigma\hat{|}Y} \colon \Omega^{p-k}_{X \hat{|} Y} \longrightarrow {\rm Gr}^W_k(\omega^p_{X\hat{|} Y})
\end{equation}
which is the completion along $Y$ of \eqref{rho}, and is locally defined as in \eqref{res05}, but with $\eta \in  \Omega^{p-k}_{X \hat{|} Y}$. Then, we can see that this map $\rho_{\sigma\hat{|}Y}$ factorizes into 
\begin{equation}
\label{res008}
\begin{matrix}
\Omega^{p-k}_{X \hat{|}Y} & \longrightarrow & i_{\sigma*} \Omega^{p-k}_{D_{\sigma} \hat{|} Y_{\sigma}} \otimes ({\rm orientation} \, \sigma_1... \sigma_k) \\
^{\rho_{\sigma}}{\downarrow} & {\swarrow}_{\overline{\rho}_{\sigma}}  &  \\
 {\rm Gr}^W_k(\omega^{p}_{X \hat{|} Y}) & & 
\end{matrix}
\end{equation}
which is the formal analogue of \eqref{res6}. We  conclude that there exists an isomorphism of complexes  on $X$
\begin{equation}
\label{res0008}
\overline{\rho}^{^.}_{\hat{|} Y} \colon \pi^k_* \Omega^{^.}_{D^k \hat{|} Y^k} (\varepsilon^k)[-k]
 \longrightarrow {\rm Gr}^W_k(\omega^{^.}_{X \hat{|} Y})
\end{equation}
whose inverse isomorphism is exactly $\hat{{\rm Res}}^{^.} $.
\end{remark}

\subsection{Cohomology of $\omega^{^.}_{X\hat{|}Y}$}
\label{res9}

From the isomorphism $\hat{{\rm Res}}^{^.} $, and from \eqref{res10}, we deduce that, when $X$ is smooth over $S$, 
with log structure given by a normal crossing divisor $D$ on $X$, and  $Y$ is a closed subscheme of $X$, 
then
\begin{equation}
\label{res09}
\mathscr H^q({\rm Gr}^W_k(\omega^{^.}_{X\hat{|}Y}))\cong \mathbb C_{Y^k} \otimes_{\mathbb C} \varepsilon^k_X 
\cong \mathbb C_{Y^k}  \otimes_{\mathbb Z} \bigwedge^k M^{gp}_X/\mathcal O_X^* \quad {\rm if} \quad q=k
\end{equation}
and 
\begin{equation}
\label{res009}
\mathscr H^q({\rm Gr}^W_k(\omega^{^.}_{X\hat{|}Y}))=0 \quad {\rm if} \quad q\neq k
\end{equation}
Therefore, we deduce that, for each point $x \in Y \cap D$, there exists an isomorphism
\begin{equation}
\label{res11}
\mathscr H^q( \omega^{^.}_{X\hat{|}Y})_x \cong \mathbb C \otimes_{\mathbb Z} \bigwedge^q (M^{gp}_X/\mathcal O_X^*)_x
\end{equation}

\section{Formal Log Poincar\'e Lemma}
\label{poin}

In this section, we generalize  the logarithmic version of the Poincar\'e Lemma, proved by Kato-Nakayama (\cite[Theorem (3.8)]{Nakayama}) in the case of an (ideally) log smooth log analytic space (i.e. a log analytic space satisfying the assumption $(0.4)$ in \cite{Nakayama}). We extend this result to the case of a  generic fs log analytic space over $S$,  and prove  the following 
\begin{theorem}
\label{poin0}
Let $i \colon Y \hookrightarrow X$ be like in \S \ref{o}, and let $\omega^{^.,log}_{X \hat{|}Y}$ be the complex introduced in Definition \ref{re9}, with differential maps given by  \eqref{re11}. Then, there exists a quasi-isomorphism 
\begin{equation}
\label{poin1}
\mathbb C_{Y^{log}} \stackrel{\cong}{\longrightarrow} \omega^{^.,log}_{X \hat{|}Y}
\end{equation}
\end{theorem}

To prove this theorem we first need some preliminary results. The methods of the proof are similar to those in \cite{Nakayama}. \\
Let $Y$ and $X$ be as in \S \ref{o}. Let $P \longrightarrow M_X$ be a chart, with $P$ a toric (or simply fs) monoid. 
Let $\mathfrak p$ be a prime ideal of $P$ which is sent to $0 \in \mathcal O_X$ under $P \longrightarrow M_X \longrightarrow \mathcal O_X$.  
Let $T$ be the fs log analytic space whose underlying space is the same as that of $X$ but whose log structure $M_T$ 
is associated to $P \smallsetminus \mathfrak p \longrightarrow \mathcal O_T$. Similarly, let $Z$ be the closed log subspace of $T$ 
whose underlying space is the same as that of $Y$ and whose log structure is  the inverse image of $M_T$. We have the following commutative diagram of fine log analytic spaces
\begin{equation}
\label{poin0005}
\begin{matrix}
(Y, i^*M_X) & \stackrel{i}{\hookrightarrow} & (X, M_X)\\
\downarrow & & \downarrow \\
(Z, i_T^* M_T) &  \stackrel{i_T}{\hookrightarrow} & (T, M_T)
\end{matrix}
\end{equation}
where the vertical maps are the identity over the underlying analytic spaces.  We also note that, since the closed immersions $i$ and $i_T$ are both exact, the log formal analytic space  $T \hat{|} Z$ coincides with the classical formal  analytic space $X \hat{|}Y$ and so 
\begin{equation}
\label{poin00005}
\omega^{^.}_{T \hat{|}Z} \cong \omega^{^.}_T \otimes_{\mathcal O_X} \mathcal O_{X \hat{|}Y}
\end{equation}
We introduce now a filtration on the complex $\omega^{^.}_{X \hat{|}Y}$. So, for $q,r \in \mathbb Z$, let $F^{\mathfrak p}_r \omega^q_X$ be the $\mathcal O_X-$subsheaf of $\omega^q_X$ defined by 
$F^{\mathfrak p}_r \omega^q_X=0$, if $r<0$; $F^{\mathfrak p}_r \omega^q_X= {\rm Im} \{ \omega^r_X \otimes \omega^{q-r}_T \longrightarrow \omega^q_X \}$, if $0 \leq r \leq q$; $F^{\mathfrak p}_r \omega^q_X=\omega^q_X$, if $q \leq r$ (\cite[${\rm Fil}_r$ in Lemma (4.4)]{Nakayama}). \\
 On the complex $\omega^{^.}_{X \hat{|}Y}$ we consider the induced filtration 
$$
\hat{F}^{\mathfrak p}_r \omega^{^.}_{X \hat{|}Y}=F^{\mathfrak p}_r \omega^{^.}_X \otimes_{\mathcal O_X} \mathcal O_{X \hat{|}Y} 
$$
\begin{Lemma}
\label{poin05}\cite[Lemma (4.4)]{Nakayama}
In the previous context, \\
(1) $\hat{F}^{\mathfrak p}_r \omega^{^.}_{X \hat{|}Y} $ is a subcomplex of $\omega^{^.}_{X \hat{|}Y}$.\\
(2) For any $r \in \mathbb Z$, there is an isomorphism of complexes
\begin{equation}
\label{poin005}
\bigwedge^r (P^{gp}/(P \smallsetminus \mathfrak p)^{gp}) \otimes_{\mathbb Z} \omega^{^.}_{T\hat{|}Z} [-r] \stackrel{\cong}{\longrightarrow} \hat{F}^{\mathfrak p}_r \omega^{^.}_{X \hat{|}Y} /\hat{F}^{\mathfrak p}_{r-1} \omega^{^.}_{X \hat{|}Y} 
\end{equation}
whose degree $q$ part is given by 
$$
(p_1 \wedge ... \wedge p_r) \otimes_{\mathbb Z} (\eta \otimes_{\mathcal O_X} f) \longmapsto {\rm dlog}(p_1) \wedge... \wedge {\rm dlog} (p_r) \wedge (\eta \otimes_{\mathcal O_X} f)
$$
where $p_1,...,p_r \in P^{gp}$, $\eta \otimes_{\mathcal O_X} f \in \omega^q_{T \hat{|}Z} \cong \omega^{^.}_{T}  \otimes_{\mathcal O_X} \mathcal O_{X \hat{|}Y} $. The differential of the left side is equal to $(id \otimes_{\mathbb Z} \hat{d}^{^.})$, where $\hat{d}^{^.} $ is the differential of $\omega^{^.}_{T \hat{|}Z}$. 
\end{Lemma}
\begin{Proof}
(2). By applying the functor $(-) \otimes \mathcal O_{X \hat{|}Y}$  to the exact sequence of coherent sheaves 4.4.1 in \cite[Lemma (4.4)]{Nakayama}, and using \eqref{poin00005}, we have 
$$
0 \longrightarrow \omega^1_{T \hat{|}Z} \longrightarrow \omega^1_{X \hat{|}Y} \longrightarrow \mathcal O_{X \hat{|}Y} \otimes_{\mathbb Z} P^{gp}/(P \smallsetminus \mathfrak p)^{gp} \longrightarrow 0
$$

\end{Proof}

Let $P$ be an fs (or toric) monoid and let $X$ be the log analytic space $ \Spec \, \mathbb C[P]$, endowed with log structure $P \longrightarrow \mathcal O_X$. Let $i \colon Y \hookrightarrow X$ be an exact closed immersion, where $Y$ is a fine log analytic space endowed with the induced log structure. We fix a point $x \in Y$. Since $i \colon Y \hookrightarrow X$ is exact, via the canonical isomorphism
$$
\omega^1_{X \hat{|}Y} \cong \mathbb C[P]_{\hat{|}Y} \otimes_{\mathbb Z} P^{gp}
$$
the map $P^{gp} \longrightarrow \omega^1_{X \hat{|}Y} $, sending $p \in P^{gp}$ to ${\rm dlog} \, p$, corresponds to the map sending $p$ to $1 \otimes p$ (\cite[\S 3]{Ogus0}). The image of this map is contained in the closed $1-$forms. Therefore, we get a map 
$$
M^{gp}_{X,x} / \mathcal O^*_{X,x} \cong P^{gp} \longrightarrow  \mathscr H^1(\omega^{^.}_{X \hat{|}Y} )
$$
and, by cup product, we deduce a map 
\begin{equation}
\label{poi9}
\bigwedge^q (M^{gp}_{X,x} / \mathcal O^*_{X,x}) \cong \bigwedge^q  P^{gp} \longrightarrow  \mathscr H^q(\omega^{^.}_{X \hat{|}Y} )
\end{equation}

Let $\mathfrak b$ be the prime ideal of $P$ which is the inverse image of the maximal ideal of $\mathcal O_{X,x}$
 under $P \longrightarrow \mathcal O_{X,x}$. 
We denote by $X(\mathfrak b)$ the closed analytic subspace $\Spec \, (\mathbb C[P]/(\mathfrak b))$ of $X$. The underlying analytic 
space of $X(\mathfrak b)$ is equal to $\Spec \,  \mathbb C[P \smallsetminus \mathfrak b]$, and $x$ belongs to its smooth open analytic 
subspace $\Spec \, \mathbb C[(P \smallsetminus \mathfrak b)^{gp}]$, where the log structure is trivial. 
Let $Y(\mathfrak b) $ be the fiber product 
\begin{equation}
\label{p03}
\begin{matrix}
Y(\mathfrak b)  & \hookrightarrow & X(\mathfrak b) \\
\downarrow & & \downarrow \\
Y & \hookrightarrow & X
\end{matrix}
\end{equation}
which is a closed subspace of $X(\mathfrak b)$. Moreover, let $(X \hat{|}Y )(\mathfrak b)$ be the completion $X(\mathfrak b) \hat{|} Y(\mathfrak b)$ of $X(\mathfrak b)$ along its closed subspace $Y(\mathfrak b)$, and let $\omega^{^.}_{(X \hat{|}Y) (\mathfrak b)}$ be the formal log complex of $(X \hat{|}Y) (\mathfrak b)$, with differential maps $\hat{d}^{^.}_{\mathfrak b} $. \\ 
We denote by  $(\mathcal O_{(X \hat{|}Y) (\mathfrak b)})^{\hat{d}^{1}_{\mathfrak b} =0}$ the kernel of $\hat{d}^1_{\mathfrak b}  \colon \mathcal O_{ (X \hat{|}Y) (\mathfrak b)} \longrightarrow \omega^{1}_{(X \hat{|}Y) (\mathfrak b)} $.  
\begin{Lemma}
\label{poin2}\cite[Lemma (4.5)]{Nakayama}
In the previous context,  \\
(1) if we restrict to some open neighbourhood of $x $ in $X(\mathfrak b)$ ($x \in Y$), there exists an isomorphism   
\begin{equation}
\label{poin3}
\mathbb C \stackrel{\cong}{\longrightarrow} (\mathcal O_{(X \hat{|}Y) (\mathfrak b)})^{\hat{d}^{1}_{\mathfrak b} =0}_x
\end{equation}
and a quasi-isomorphism 
\begin{equation}
\label{poin03}
\mathbb C \stackrel{\cong}{\longrightarrow} (\omega^{^.}_{(X \hat{|}Y) (\mathfrak b)})_x
\end{equation}
(2) Let $\mathfrak b \omega^{^.}_{X \hat{|}Y}$ be the subcomplex of $\omega^{^.}_{X \hat{|}Y}$ whose degree $q$ part is defined to be the $\mathcal O_X-$subsheaf of $\omega^{q}_{X \hat{|}Y}$ generated by $ b \omega^{q}_{X \hat{|}Y}$, with $b \in \mathfrak b$. For any $q$, the map 
\begin{equation}
\label{poin4}
\bigwedge^q (M_X^{gp}/\mathcal O_X^*)_x \otimes_{\mathbb Z} \mathbb C \longrightarrow \mathscr H^q(\omega^{^.}_{X \hat{|}Y}/ \mathfrak b \omega^{^.}_{X \hat{|}Y})_x
\end{equation}
which is induced by the map \eqref{poi9}, is bijective.\\
(3) The stalk at $x$ of the canonical map of complexes 
\begin{equation}
\label{poin5}
\omega^{^.}_{X \hat{|}Y} \longrightarrow \omega^{^.}_{X \hat{|}Y}/ \mathfrak b \omega^{^.}_{X \hat{|}Y}
\end{equation}
is a quasi-isomorphism. 
\end{Lemma}
\begin{Proof}
We first note that the complex $\omega^{^.}_{X \hat{|}Y}/ \mathfrak b \omega^{^.}_{X \hat{|}Y}$ is isomorphic to $\omega^{^.}_{(X \hat{|}Y) (\mathfrak b)}$. Indeed, since $\mathfrak b$ is an ideal of the monoid $P$, by \cite[Lemma (3.6), (2)]{Nakayama}, $\omega^{1}_X/ \mathfrak b \omega^{1}_X \cong \omega^1_{X(\mathfrak b)}$, and it follows that 
$$
\omega^{^.}_{X \hat{|}Y}/ \mathfrak b \omega^{^.}_{X \hat{|}Y} \cong \omega^{^.}_{X}/ \mathfrak b \omega^{^.}_{X} \otimes_{\mathcal O_X} \mathcal O_{X \hat{|}Y} \cong 
$$
$$
 \omega^{^.}_{X(\mathfrak b)} \otimes_{\mathcal O_X} \mathcal O_{X \hat{|}Y} \cong  \omega^{^.}_{X(\mathfrak b) \hat{|} Y(\mathfrak b)}= \omega^{^.}_{(X \hat{|}Y) (\mathfrak b)}
$$
We start to prove (1) and (2). We may restrict ourselves to the  open neighbourhood 
$\Spec \, \mathbb C[(P \smallsetminus \mathfrak b)^{gp}]$ of $x$ in $X(\mathfrak b)$, 
and consider the restriction of $Y(\mathfrak b)$ to this open neighborhood. So, $x$  belongs to the 
closed subspace $Y(\mathfrak b)\cap \Spec \, \mathbb C[(P \smallsetminus \mathfrak b)^{gp}]$ 
of the non-singular analytic space $\Spec \, \mathbb C[(P \smallsetminus \mathfrak b)^{gp}]$, 
where the log structure is trivial. \\
In this local situation,  $\omega^{^.}_{(X \hat{|}Y) (\mathfrak b)} \cong \Omega^{^.}_{X(\mathfrak b) \hat{|} Y(\mathfrak b)}$, 
and, from \cite[Theorem 2.1]{Hartshorne}, we  know that this complex is a resolution of 
the constant sheaf $\mathbb C_{Y(\mathfrak b)}$ over $Y(\mathfrak b)$. Therefore, the stalk at $x $ 
of $(\mathcal O_{(X \hat{|}Y) (\mathfrak b)})^{\hat{d}^{1}_{\mathfrak b} =0}$ is isomorphic to $\mathbb C$, and there is a quasi-isomorphism
$$
\mathbb C \stackrel{\cong}{\longrightarrow} (\Omega^{^.}_{X(\mathfrak b) \hat{|} Y(\mathfrak b)})_x
$$
so (1) is proved. Now, we apply Lemma \ref{poin05} by taking  $X(\mathfrak b)$, $P$, $\mathfrak b$ as $X$, $P$ and $\mathfrak p$. 
We consider $\mathscr H^q$ of both sides of Lemma \ref{poin05}, (2), and take the stalk at $x$. Then, 
$\mathscr H^q (\hat{F}^{\mathfrak b}_r(\omega^{^.}_{(X \hat{|}Y)(\mathfrak b)})/\hat{F}^{\mathfrak b}_{r-1}(\omega^{^.}_{(X \hat{|}Y)(\mathfrak b)}))_x$
 is isomorphic to 
$\bigwedge^r P^{gp}/(P \smallsetminus \mathfrak b)^{gp} \otimes_{\mathbb Z} \mathscr H^{q-r}(\Omega^{^.}_{X(\mathfrak b) \hat{|} Y(\mathfrak b)})_x $,
 which is isomorphic to   
 $\bigwedge^r P^{gp}/(P \smallsetminus \mathfrak b)^{gp} \otimes_{\mathbb Z} \mathbb C$ if $q=r$ 
and is zero if $q \neq r$. 
Therefore, since $\omega^{^.}_{X \hat{|}Y}/ \mathfrak b \omega^{^.}_{X \hat{|}Y} \cong \omega^{^.}_{(X \hat{|}Y) (\mathfrak b)}$, 
the stalk at $x$ of $\mathscr H^q(\omega^{^.}_{X \hat{|}Y}/ \mathfrak b \omega^{^.}_{X \hat{|}Y})$ is isomorphic 
to $\bigwedge^q (M^{gp}_{X,x}/\mathcal O^*_{X,x}) \otimes_{\mathbb Z} \mathbb C$, so (2) is proved.\\
To prove (3), we consider the particular case where $P= \mathbb N^r$, for some $r \geq 0$. In this situation, $X \cong \mathbb C^r$ as an 
analytic  space, with canonical log structure given by a normal crossing divisor $D \hookrightarrow X$, and $Y $ is a closed analytic subspace of $X$,
 with induced log structure. 
Then, the complex $\omega^{^.}_{X\hat{|}Y} \cong \Omega^{^.}_X({\rm log} \, D) \otimes_{\mathcal O_X} \mathcal O_{X \hat{|}Y}$. 
Therefore, we are reduced to the case analyzed in \S \ref{res9}, and we can use the isomorphism \eqref{res11} 
to describe the stalk at $x $ of $\mathscr H^q(\omega^{^.}_{X\hat{|}Y} )$. So,  by applying Lemma \ref{poin2}, (2), 
$\mathscr H^q(\omega^{^.}_{X \hat{|}Y}/ \mathfrak b \omega^{^.}_{X \hat{|}Y})_x \cong 
 \bigwedge^q (M_X^{gp}/\mathcal O_X^*)_x \otimes_{\mathbb Z} \mathbb C$, which is isomorphic to $\mathscr H^q(\omega^{^.}_{X\hat{|}Y} )_x$, via \eqref{res11}. \\
Now, we prove (3) in the general situation. For a non-empty ideal $I$ of the monoid $P$, we can consider 
the toric variety $B_I(\Spec \, \mathbb C[P])$, which we get from $X$ by ``blowing-up" along $I$ as in \cite[\S I, Theorem 10]{Kempf}.
It is endowed with a canonical log structure (\cite[(3.7)(1)]{Kato1}). \\

{\bf Note.}  From \cite[Proposition (9.8)]{Kato2}, and \cite[\S I, Theorem 11]{Kempf}, it is possible to choose an ideal $\tilde{I}$ of $P$, 
such that, if $\tilde{X}= B_{\tilde{I}}(\Spec \, \mathbb C[P])$, with log structure $\tilde{M}$, 
then, for any $y \in \tilde{X}$, $(\tilde{M}/\mathcal O^*_{\tilde{X}})_y$ is isomorphic to $\mathbb N^{r(y)}$,
 for some $r(y) \geq 0$. Let $f \colon \tilde{X} \longrightarrow X$ be the  proper  map, corresponding to the ``blowing-up" of $X$ 
along $\tilde{I}$. We can suppose $\tilde{X} \cong \mathbb C^r$, for some $r \in \mathbb N$, i.e. $X$ to be a
 non-singular log analytic space, with canonical log structure $\mathbb N^r \longrightarrow \mathbb C[\mathbb N^r]$,
namely, the log structure given by a normal crossing divisor $\tilde{D}$. \\

Then, we consider the following cartesian diagram 
\begin{equation}
\label{poi0}
\begin{matrix}
\tilde{Y } & \stackrel{\tilde{i}}{\hookrightarrow}  & \tilde{X}\\
^{f_{|Y}}{\downarrow} & & ^{f}{\downarrow} \\
Y & \stackrel{i}{\hookrightarrow} & X
\end{matrix}
\end{equation}
where $\tilde{Y}=f^{-1}(Y)$ is a closed subspace of $\tilde{X}$, and we suppose it to be endowed with the inverse image of the 
log structure $\tilde{M}$. We denote by $\hat{f}$ the morphism 
$\hat{f} \colon \tilde{X} \hat{|} \tilde{Y} \longrightarrow X \hat{|} Y$ (deduced from the cartesian diagram \eqref{poi0}). 
We also note that the vertical maps in \eqref{poi0} are log-\'etale, so 
\begin{equation}
\label{poi1}
\omega^{^.}_{\tilde{X}} \cong f^* \omega^{^.}_X
\end{equation}
Then, from \eqref{poi1},  we get 
\begin{equation}
\label{poi01}
\omega^{^.}_{\tilde{X} \hat{|} \tilde{Y}} \cong \hat{f}^* \omega^{^.}_{X\hat{|}Y} 
\end{equation}
Moreover, by \cite[\S I, Corollary 1. c)]{Kempf}, there exists a  quasi-isomorphism
\begin{equation}
\label{poi2}
\mathcal O_X \stackrel{\cong}{\longrightarrow } \mathbb Rf_* \mathcal O_{\tilde{X}}
\end{equation} 
Since $f$ is proper, and $X$, $\tilde{X} $ are schemes of finite type over $S$, 
applying \cite[\S II, Proposition (6.2)]{Hartshorne} to the structural sheaf $\mathcal O_{\tilde{X}}$, we get 
\begin{equation}
\label{poi3}
(\mathbb R f_* \mathcal O_{\tilde{X}})_{\hat{|}Y} \stackrel{\cong}{\longrightarrow} \mathbb R \hat{f}_* (\mathcal O_{\tilde{X} \hat{|} \tilde{Y}})
\end{equation}
and, from the isomorphism \eqref{poi2}, we get 
\begin{equation}
\label{poi4}
\mathcal O_{X \hat{|} Y} \stackrel{\cong}{\longrightarrow} \mathbb R \hat{f}_* (\mathcal O_{\tilde{X} \hat{|} \tilde{Y}})
\end{equation}
Therefore, since the $\mathcal O_{\tilde{X}}$-module (resp. $\mathcal O_{X }-$module) $\omega^q_{\tilde{X}}$ (resp. $\omega^q_X$)
 is free of finite rank, for any $q$, from \eqref{poi01} and \eqref{poi4}, we finally get an isomorphism in the derived category
\begin{equation}
\label{poi5}
\omega^{^.}_{X \hat{|}Y} \cong \mathbb R \hat{f}_* \omega^{^.}_{\tilde{X} \hat{|} \tilde{Y}}
\end{equation}
Let now $\tilde{x} \in \tilde{Y}$ be such that $f(\tilde{x})=x$. Let $\tilde{\mathfrak b}$ be the prime ideal of $\mathbb N^r$ equal to the inverse
 image of the maximal ideal of $\mathcal O_{\tilde{X},x}$ under $\mathbb N^r \longrightarrow \mathcal O_{\tilde{X}, x}$.
 We consider the log analytic closed subspace $\tilde{X}(\tilde{\mathfrak b})$ of $\tilde{X}$,
 and its closed log subspace $\tilde{Y} (\tilde{\mathfrak b}) \hookrightarrow \tilde{X}(\tilde{\mathfrak b})$, defined as in \eqref{p03}. \\
We consider the following commutative diagram 
\begin{equation}
\label{poi6}
\begin{matrix}
\omega^{^.}_{X \hat{|}Y} & \longrightarrow & \omega^{^.}_{X \hat{|}Y}/ \mathfrak b \omega^{^.}_{X \hat{|}Y}
 \cong \omega^{^.}_{(X \hat{|}Y)(\mathfrak b)} \\
^{\cong}{\downarrow} & & \downarrow \\
\mathbb R \hat{f}_* \omega^{^.}_{\tilde{X} \hat{|} \tilde{Y}} & \longrightarrow
 & \mathbb R \hat{f}_*( \omega^{^.}_{\tilde{X} \hat{|}\tilde{Y}}/ \tilde{ \mathfrak b} \omega^{^.}_{\tilde{X} \hat{|} \tilde{Y}}) 
\cong \mathbb R \hat{f}_*  
\omega^{^.}_{(\tilde{X} \hat{|} \tilde{Y})(\tilde{\mathfrak b})}
\end{matrix}
\end{equation}
Since we have just proved (3) in the case $P= \mathbb N^r$, $r \geq 0$, the lower horizontal arrow is an isomorphism at $x$. Therefore, the map 
\begin{equation}
\label{p00}
\mathscr H^q (\omega^{^.}_{X \hat{|}Y})_x \longrightarrow \mathscr H^q( \omega^{^.}_{X \hat{|}Y}/ \mathfrak b \omega^{^.}_{X \hat{|}Y} )_x
\end{equation}
is injective. Moreover, by Lemma \ref{poin2}, (1), 
$\mathscr H^0(\omega^{^.}_{X \hat{|}Y}/ \mathfrak b \omega^{^.}_{X \hat{|}Y})_x \cong
 \mathscr H^0(\omega^{^.}_{(X \hat{|}Y)(\mathfrak b)})_x \stackrel{\cong}{\longrightarrow} \mathbb C$, 
and also $\mathscr H^0(\omega^{^.}_{\tilde{X} \hat{|}\tilde{Y}}/ \tilde{ \mathfrak b} \omega^{^.}_{\tilde{X} \hat{|} \tilde{Y}} )
 \cong \mathscr H^0(\omega^{^.}_{(\tilde{X} \hat{|} \tilde{Y})(\tilde{\mathfrak b})}) \cong \mathbb C$. 
Thus, we find that the composition map 
$$
\mathbb C \cong \mathscr H^0(\omega^{^.}_{X \hat{|}Y}/ \mathfrak b \omega^{^.}_{X \hat{|}Y})_x \longrightarrow \mathscr H^0(\mathbb R \hat{f}_* (\omega^{^.}_{\tilde{X} \hat{|}\tilde{Y}}/ \tilde{ \mathfrak b} \omega^{^.}_{\tilde{X} \hat{|} \tilde{Y}} ))_x 
\longrightarrow \mathscr H^0(\omega^{^.}_{\tilde{X} \hat{|}\tilde{Y}}/ \tilde{ \mathfrak b} \omega^{^.}_{\tilde{X} \hat{|} \tilde{Y}} )_{\tilde{x}} \cong \mathbb C
$$
is the identity map, and so 
\begin{equation}
\label{poi7}
\mathscr H^0(\omega^{^.}_{X \hat{|}Y}/ \mathfrak b \omega^{^.}_{X \hat{|}Y})_x
 \longrightarrow
 \mathscr H^0(\mathbb R \hat{f}_* (\omega^{^.}_{\tilde{X} \hat{|}\tilde{Y}}/ \tilde{ \mathfrak b} \omega^{^.}_{\tilde{X} \hat{|} \tilde{Y}} ))_x 
\end{equation}
is injective. Now, from diagram \eqref{poi6}, since the composed map 
$\omega^{^.}_{X \hat{|}Y} \longrightarrow
 \mathbb R \hat{f}_*( \omega^{^.}_{\tilde{X} \hat{|}\tilde{Y}}/ \tilde{ \mathfrak b} \omega^{^.}_{\tilde{X} \hat{|} \tilde{Y}}) $ 
is an isomorphism at $x$, it follows that the map \eqref{poi7} is also surjective, and so it is an isomorphism.\\
Therefore, from diagram \eqref{poi6}, the map   
\begin{equation}
\label{poi8}
\mathscr H^0(\omega^{^.}_{X \hat{|}Y})_x \longrightarrow \mathscr H^0(\omega^{^.}_{X \hat{|}Y}/ \mathfrak b \omega^{^.}_{X \hat{|}Y})_x
\end{equation}
is also an isomorphism, and we can conclude that  $\mathscr H^0(\omega^{^.}_{X \hat{|}Y})_x  \cong \mathbb C$. \\
Now, the isomorphism \eqref{poin4}, factorizes through 
$$
\bigwedge^q (M^{gp}_{X,x}/\mathcal O^*_{X,x})\otimes_{\mathbb Z} \mathscr H^0(\omega^{^.}_{X \hat{|}Y})_x
 \longrightarrow \mathscr H^q(\omega^{^.}_{X \hat{|}Y})_x \longrightarrow 
 \mathscr H^q (\omega^{^.}_{X \hat{|}Y}/ \mathfrak b \omega^{^.}_{X \hat{|}Y})_x
$$
and we can conclude that the map \eqref{p00} is also surjective, and so it is an isomorphism.
\end{Proof}

From Lemma \ref{poin2}, we can deduce the following
\begin{proposition}
\label{poin6}\cite[Proposition (4.6)]{Nakayama}
Let $Y$ be an fs log analytic space over $S$, and let $i \colon Y \hookrightarrow X$ be an exact closed immersion 
of $Y$ into an fs log smooth log analytic space $X$. Then, for all $q \in \mathbb Z$, there is an isomorphism
\begin{equation}
\label{poin06}
\bigwedge^q (M_X^{gp}/\mathcal O_X^*)_{|Y} \otimes_{\mathbb Z} \mathbb C \stackrel{\cong}{\longrightarrow} \mathscr H^q(\omega^{^.}_{X \hat{|}Y})
\end{equation}
induced by the map $ \hat{{\rm dlog}} \colon M^{gp}_{X\hat{|}Y} \longrightarrow  \omega^1_{X \hat{|}Y}$.
\end{proposition}
\begin{Proof}
Since the question is local on $X$, we may assume that $X= \Spec \, \mathbb C[P]$, where $P $ is an fs monoid. 
Let $x \in Y$, and let $\mathfrak b \subset P$ be the inverse image of the maximal ideal of $\mathcal O_{X,x}$. Now, by Lemma \ref{poin2}, (3), 
$$
\mathscr H^q(\omega^{^.}_{X \hat{|}Y})_x \cong \mathscr H^q (\omega^{^.}_{X \hat{|}Y}/ \mathfrak b \omega^{^.}_{X \hat{|}Y})_x
$$
and, by Lemma \ref{poin2}, (2),
$$
\mathscr H^q (\omega^{^.}_{X \hat{|}Y}/ \mathfrak b \omega^{^.}_{X \hat{|}Y})_x \cong \bigwedge^q (M_{X,x}^{gp}/\mathcal O_{X,x}^*) \otimes_{\mathbb Z} \mathbb C
$$
for each point $x \in Y$.
\end{Proof}

Now, we use Lemma \ref{poin2} and Proposition \ref{poin6} to prove a ``formal version" of the logarithmic Poincar\'e Lemma.

\begin{Proof} (of Theorem \ref{poin0}) \\
In the previous notation, let $x \in Y$, $y \in Y^{log}$ be such that $\tau(y)=x$. Let $\{ t_1,...,t_n\} $ be a family of
 elements of $\mathscr L_{X,x}$ whose image via the map ${\rm exp}_x \colon \mathscr L_{X,x} \longrightarrow M_{X,x}^{gp}/\mathcal O_{X,x}^*$ 
is a $\mathbb Z-$basis of $M_{X,x}^{gp}/\mathcal O_{X,x}^*$. \\
Let $R$ be the polynomial ring $ \mathbb C[T_1,...,T_n] $. From Lemma \ref{re6}, the stalk in $y$ of $\mathcal O^{log}_{X \hat{|}Y}$ 
is isomorphic to $\mathcal O_{X \hat{|}Y,x} [T_1,...,T_n]$, where each $t_i$ corresponds to $T_i$. Therefore, we consider
 the $\mathbb C-$linear homomorphism 
\begin{equation}
\label{poin7}
R \longrightarrow \mathcal O^{log}_{X \hat{|}Y,y}
\end{equation}
which sends $T_i \longmapsto t_i$,  for $i=1,...,n$. Since 
$$
\mathbb C \longrightarrow \Omega^{^.}_{R/\mathbb C}
$$
is a quasi-isomorphism, it is sufficient to prove that there exists a quasi-isomorphism
\begin{equation}
\label{poin8}
\Omega^{^.}_{R/\mathbb C} \stackrel{\cong}{\longrightarrow} \omega^{^.,log}_{X \hat{|}Y,y}
\end{equation}
To this end, we introduce a filtration on $\Omega^{^.}_{R /\mathbb C}$ as follows:  for any $r \in \mathbb Z$, 
let ${\rm Fil}_r(\Omega^{^.}_{R /\mathbb C})$ be the subcomplex of $\Omega^{^.}_{R /\mathbb C}$ whose degree $q$ part 
is the $\mathbb C-$submodule of $\Omega^{q}_{R /\mathbb C}$ 
generated by elements of the type $f \cdot \gamma$, with $f \in R$ an element of degree $\leq r$, and $\gamma \in \Omega^{q}_{R /\mathbb C}$. \\
We also introduce a filtration on $\omega^{^.,log}_{X \hat{|}Y,y}$: for any $r \in \mathbb Z$, 
let ${\rm Fil}_r(\omega^{^.,log}_{X \hat{|}Y})$ be the subcomplex of $\omega^{^.,log}_{X \hat{|}Y}$ 
whose degree $q$ part ${\rm Fil}_r(\omega^{q,log}_{X \hat{|}Y})$ is defined as 
$$
{\rm Fil}_r(\omega^{q,log}_{X \hat{|}Y})=: \hat{{\rm fil}}_r(\mathcal O^{log}_{X \hat{|}Y}) \otimes_{\tau^{-1}(\mathcal O_X)} \tau^{-1}(\omega^q_X)
$$
where $\hat{{\rm fil}}_r(\mathcal O^{log}_{X \hat{|}Y})$ is the filtration defined in Lemma \ref{re8}. \\
Then, by Lemma \ref{re8},
$$
{\rm Fil}_r(\omega^{^.,log}_{X \hat{|}Y})/{\rm Fil}_{r-1}(\omega^{^.,log}_{X \hat{|}Y}) \cong \tau^{-1}\left( \omega^{^.}_{X \hat{|}Y}
 \otimes_{\mathbb Z} {\rm Sym}^r_{\mathbb Z}(M_X^{gp}/\mathcal O_X^*) \right)
$$
and, by Proposition \ref{poin6}, for any $q$,
$$
\mathscr H^q \left( \tau^{-1} ( \omega^{^.}_{X \hat{|}Y} \otimes_{\mathbb Z} {\rm Sym}^r_{\mathbb Z}(M_X^{gp}/\mathcal O_X^*)) \right) \cong 
 \mathbb C \otimes_{\mathbb Z} \bigwedge^q(M_X^{gp}/\mathcal O_X^*) \otimes_{\mathbb Z} \tau^{-1}({\rm Sym}^r_{\mathbb Z}(M_X^{gp}/\mathcal O_X^*))
$$
On the other hand, ${\rm Fil}_r(\Omega^{^.}_{R /\mathbb C})/ {\rm Fil}_{r-1}(\Omega^{^.}_{R /\mathbb C})$ is the complex  
$$
q \longmapsto \mathbb C \otimes_{\mathbb Z} (\bigwedge^q \bigoplus_{i=1}^n \mathbb Z T_i) 
\otimes_{\mathbb Z} {\rm Sym}_{\mathbb Z}^r(\bigoplus_{i=1}^n \mathbb Z T_i )
$$
which is isomorphic to $\mathbb C \otimes_{\mathbb Z} \bigwedge^q(M_X^{gp}/\mathcal O_X^*)_x \otimes_{\mathbb Z} 
{\rm Sym}^r_{\mathbb Z}(M_X^{gp}/\mathcal O_X^*)_x$. The differentials of this complex are zero. \\
Therefore, for any $r \in \mathbb Z$, there exists a quasi-isomorphism
\begin{equation}
\label{poin9}
{\rm Fil}_r(\Omega^{^.}_{R /\mathbb C})/ {\rm Fil}_{r-1}(\Omega^{^.}_{R /\mathbb C}) \stackrel{\cong}{\longrightarrow}  
{\rm Fil}_r(\omega^{^.,log}_{X \hat{|}Y})/{\rm Fil}_{r-1}(\omega^{^.,log}_{X \hat{|}Y}) 
\end{equation}
and this implies that there is a quasi-isomorphism $\Omega^{^.}_{R /\mathbb C} \cong \omega^{^.,log}_{X \hat{|}Y,y}$ for each point $y \in Y^{log}$.
\end{Proof}

\section{Log De Rham and Log Betti Cohomologies}
\label{finale}

The goal of this section  is to compare the Log Betti Cohomology 
$H^{^.}(Y^{log}, \mathbb C)$ of an fs log scheme $Y$, with its Algebraic Log  De Rham Cohomology $\mathbb H^{^.}(Y, \omega^{^.}_{X \hat{|} Y})$.
 Therefore, we begin with
 \begin{theorem}
\label{f00}
Let $Y$ be an fs log scheme over $S$, and let $i \colon Y \hookrightarrow X$ be an exact closed immersion of $Y$ 
into an fs log smooth log scheme $X$.  
 Then, for any $q \in \mathbb Z$, there exists  an isomorphism
\begin{equation}
\label{f01}
H^q(Y^{log}, \mathbb C) \cong \mathbb H^q(Y, \omega^{^.}_{X \hat{|} Y})
\end{equation}
\end{theorem}
\begin{Proof}
In the previous section, we have checked that 
$H^q(Y^{log}, \mathbb C) \cong \mathbb H^q(Y^{log}, \omega^{^.,log}_{X \hat{|}Y})$, for any $q \in \mathbb Z$. So,  we will first show that the Log Betti Cohomology of $Y$ is isomorphic to the Analytic Log De Rham Cohomology $\mathbb H^{^.} (Y^{an}, \omega^{^.}_{(X \hat{|}Y)^{an}})$ (Proposition \ref{f0}). Finally, we will check that the algebraic log  De Rham complex $\omega^{^.}_{X \hat{|} Y}$ is quasi-isomorphic to its associated analytic log complex $\omega^{^.}_{(X \hat{|}Y)^{an}}$ (Theorem  \ref{f3}).
\end{Proof} 
\begin{proposition}
\label{f0} \cite[(4.8), 4.8.5]{Nakayama}. There exists a quasi-isomorphism
\begin{equation}
\label{f1}
\mathbb R \tau_* (\omega^{^.,log}_{X \hat{|}Y}) \stackrel{\cong}{\longrightarrow} \omega^{^.}_{(X \hat{|}Y)^{an}}
\end{equation}
\end{proposition}
\begin{Proof}
We consider the composed map
\begin{equation}
\label{f2}
\mathbb C \otimes_{\mathbb Z} \bigwedge^q (M_X^{gp}/ \mathcal O_X^*)_{|Y} \stackrel{\cong}{\longrightarrow} 
\mathscr H^q(\omega^{^.}_{(X \hat{|}Y)^{an}}) \longrightarrow \mathbb R^q \tau_* \mathbb C_{Y^{log}} 
\end{equation}
(where the second map is given by the natural arrow $\omega^{^.}_{(X \hat{|}Y)^{an}} \longrightarrow \mathbb R \tau_* 
 (\omega^{^.,log}_{X \hat{|}Y})$). Since the first map is an isomorphism by Proposition \ref{poin6}, 
it is sufficient to prove that there is a quasi-isomorphism
$$
\mathbb R^q \tau_* \mathbb C_{Y^{log}} \cong \mathbb C \otimes_{\mathbb Z} \bigwedge^q(M_X^{gp}/ \mathcal O_X^*)_{|Y} 
$$
To prove this, we apply \cite[Lemma (1.5)]{Nakayama}, taking the constant sheaf $\mathbb C$ on $Y^{an}$.  We have a canonical isomorphism
$$
(\mathbb R^q \tau_* \mathbb C_{Y^{log}} \cong )\mathbb R^q \tau_* \tau^* \mathbb C_{Y^{an}} \cong \mathbb C \otimes_{\mathbb Z}
 \bigwedge^q M_Y^{gp}/\mathcal O_Y^*
$$
where the sheaf $M_Y^{gp}/\mathcal O_Y^*$ is isomorphic to $(M_X^{gp}/ \mathcal O_X^*)_{|Y} $, by \eqref{re1}.
\end{Proof}\\
We will now compare the algebraic log De Rham complex $\omega^{^.}_{X \hat{|}Y}$ with its associated analytic 
log De Rham complex $\omega^{^.}_{(X \hat{|}Y)^{an}}$, and we will show that they are quasi-isomorphic. 
\begin{theorem}
\label{f3}
In the previous notation, let $g \colon X^{an} \longrightarrow X$ be the canonical morphism. If we consider the cartesian diagram 
\begin{equation}
\label{f03}
\begin {matrix}
Y^{an} & \stackrel{i^{an}}{\hookrightarrow} & X^{an} \\
^{g_Y}{\downarrow} & & ^{g}{\downarrow} \\
Y & \stackrel{i}{\hookrightarrow} & X    
\end{matrix}
\end{equation}
then the morphism 
\begin{equation}
\label{f04}
\omega^{^.}_{X \hat{|}Y} \longrightarrow \mathbb R \hat{g}_* \omega^{^.}_{(X \hat{|}Y)^{an}}
\end{equation} 
induces an isomorphism in cohomology 
\begin{equation}
\label{f4}
\mathbb H^{^.}(Y, \omega^{^.}_{X \hat{|}Y}) \stackrel{\cong}{\longrightarrow} \mathbb H^{^.}(Y^{an},  \omega^{^.}_{(X \hat{|}Y)^{an}})
\end{equation}
\end{theorem}  
\begin{Proof}
We may assume that $X= \Spec \, \mathbb C[P]$, endowed with the canonical log structure (with $P$ a toric monoid), and $Y \hookrightarrow X$ 
is a closed log subscheme of $X$, with the induced log structure. We divide the proof into two steps: \\
1) We begin by proving the assertion in the case where $P= \mathbb N^r$, for some $r \in \mathbb N$, i.e. in the case of a smooth scheme 
$X$ over $S$, with log structure given by a normal crossing divisor $D \hookrightarrow X$. Then, by the formal Poincar\'e residue 
isomorphism \eqref{res5}, for each $k \leq n$, we have the following  identifications
\begin{equation}
\label{f5}
\mathbb H^q (Y,  {\rm Gr}^W_k (\omega^{^.}_{X \hat{|}Y})) \cong \mathbb H^{q-k} (Y, \hat{\pi}^k_* \Omega^{^.}_{D^k \hat{|} Y^k} (\varepsilon^k)) 
\end{equation}
Moreover, by \cite[\S IV]{Hartshorne},
$$
\mathbb H^{q-k} (Y, \hat{\pi}^k_* \Omega^{^.}_{D^k \hat{|} Y^k} (\varepsilon^k)) \cong 
\mathbb H^{q-k} (Y^{an}, \hat{\pi}^k_* \Omega^{^.}_{D^{k,an} \hat{|} Y^{k,an}} (\varepsilon^k)) \cong H^{q-k}(Y^{k,an}, \mathbb C)
$$  
and so 
$$
\mathbb H^q (Y,  {\rm Gr}^W_k (\omega^{^.}_{X \hat{|}Y})) \cong \mathbb H^q (Y^{an},  {\rm Gr}^W_k (\omega^{^.}_{(X \hat{|}Y)^{an}})) 
$$
for each $k$, $0 \leq k \leq n$. 
Therefore, we can conclude that the morphism (\ref{f04}) induces the isomorphism 
$\mathbb H^{^.}(Y, \omega^{^.}_{X \hat{|}Y}) \cong \mathbb H^{^.}(Y^{an},  \omega^{^.}_{(X \hat{|}Y)^{an}})$.\\
2) We now prove  the assertion for a generic toric monoid $P$.  We can take $\tilde{I}$, $\tilde{X}$, 
and $f \colon \tilde{X} \longrightarrow X $, 
as in the {\bf Note} interpolated in  the proof of Lemma \ref{poin2}. We consider the cartesian diagram \eqref{poi0}, 
for the algebraic and analytic cases. 
Then, by arguments similar to those  in the proof of Lemma \ref{poin2}, \eqref{poi01}, \eqref{poi4}, \eqref{poi5}, 
we can conclude that, in the algebraic setting,
$$
\omega^{^.}_{X \hat{|}Y} \cong \mathbb R \hat{f}_* \omega^{^.}_{\tilde{X} \hat{|} \tilde{Y}} 
$$
and similarly, in the analytic setting, 
$$
\omega^{^.}_{(X \hat{|}Y)^{an}} \cong \mathbb R \hat{f}^{an}_* \omega^{^.}_{(\tilde{X} \hat{|} \tilde{Y})^{an}}
$$
Therefore, to prove the assertion it is sufficient to check that there exists an isomorphism
 $\mathbb H^{^.}(\tilde{Y}, \omega^{^.}_{\tilde{X} \hat{|} \tilde{Y}}) 
\cong \mathbb H^{^.}(\tilde{Y}^{an},  \omega^{^.}_{(\tilde{X} \hat{|} \tilde{Y})^{an}})$. But this
 follows from step 1), because $\tilde{X} = \Spec \, \mathbb C[\mathbb N^r]$, for some $r \in \mathbb N$, 
endowed with canonical log structure $\mathbb N^r \longrightarrow \mathbb C[\mathbb N^r]$. 
\end{Proof}

\end{document}